\newif\ifTwoColumn
\newif\ifTechReport

\TechReporttrue    

\ifTwoColumn
\documentclass[journal]{IEEEtran}
\else
\ifTechReport
\documentclass[journal,12pt,onecolumn,draftcls]{IEEEtran}
\else
\documentclass[journal,11pt,onecolumn,draftclsnofoot]{IEEEtran}
\fi
\fi

\usepackage[utf8]{inputenc} 
\usepackage[T1]{fontenc}    
\usepackage{url}            
\usepackage{booktabs}       
\usepackage{amsfonts,amsmath,amsthm,amssymb,mathrsfs}       
\usepackage{bbm}
\usepackage{cite}
\usepackage{mathtools,xparse}
\usepackage[makeroom]{cancel}
\usepackage[acronym]{glossaries}
\usepackage{color}  
\usepackage{balance}
\usepackage[ruled,vlined]{algorithm2e}
\usepackage{algpseudocode}
\usepackage{subcaption}
\usepackage{nicefrac}       
\usepackage{microtype}      
\usepackage{xcolor}         
\usepackage{graphicx}
\usepackage[prependcaption,colorinlistoftodos]{todonotes}

\renewcommand{\qedsymbol}{$\blacksquare$}
\newcommand*{\QEDB}{\hfill\ensuremath{\square}}%

\DeclareMathSizes{10}{10}{6}{4}

\newcommand{\normaltext}[1]{\textnormal{#1}}

\graphicspath{{figures/}}

\newcommand{\R}{\mathbb{R}}

\newcommand{\N}{\mathbb{N}}

\newcommand{\bbS}{\mathbb{S}}
\newcommand{\mc}[1]{\mathcal{#1}}

\newcommand{\prob}{\mathbb{P}}

\newcommand{\col}{\mathrm{col}}

\newcommand{\bs}[1]{\boldsymbol{#1}}

\newtheorem{theorem}{Theorem}

\newtheorem{definition}{Definition}
\newtheorem{proposition}{Proposition}
\newtheorem{lemma}{Lemma}

\newtheorem{remark}{Remark}
\newtheorem{assumption}{Assumption}
\newtheorem{standing}{Standing Assumption}

\makeglossaries
\newacronym{}{}{}

\makeglossaries
\newacronym{NEP}{NEP}{Nash equilibrium problem}
\newacronym{GNEP}{GNEP}{generalized Nash equilibrium problem}
\newacronym{iid}{i.i.d.\@}{independent and identically distributed}
\newacronym{wrt}{w.r.t.\@}{with respect to}
\newacronym{psd}{psd}{positive semidefinite}
\newacronym{ls}{LS}{least squares}
\newacronym{gp}{GP}{Gaussian process}
\newacronym{maas}{MaaS}{mobility as a service}

\newglossaryentry{VI}
{
	name={VI},
	description={variational inequality},
	first={\glsentrydesc{VI} (\glsentrytext{VI})},
	plural={VIs},
	descriptionplural={variational inequalities},
	firstplural={\glsentrydescplural{VI} (VIs)}
}

\newglossaryentry{v-GNE}
{
	name={v-GNE},
	description={variational generalized Nash equilibrium},
	first={\glsentrydesc{v-GNE} (\glsentrytext{v-GNE})},
	plural={v-GNE},
	descriptionplural={variational generalized Nash equilibria},
	firstplural={\glsentrydescplural{v-GNE} (\glsentryplural{v-GNE})}
}

\newglossaryentry{GNE}
{
	name={GNE},
	description={generalized Nash equilibrium},
	first={\glsentrydesc{GNE} (\glsentrytext{GNE})},
	plural={GNE},
	descriptionplural={generalized Nash equilibria},
	firstplural={\glsentrydescplural{GNE} (\glsentryplural{GNE})}
}

\begin{document}

\title{Personalized incentives as feedback design in generalized Nash equilibrium problems}
	\author{Filippo Fabiani, Andrea Simonetto, and Paul J. Goulart
	\thanks{At the time of the work, F. Fabiani was with the Department of Engineering Science, University of Oxford, OX1 3PJ, United Kingdom.
		Currently, he is with the IMT School for Advanced Studies Lucca, Piazza San Francesco 19, 55100 Lucca, Italy ({\tt \footnotesize filippo.fabiani@imtlucca.it}). A. Simonetto is with the UMA, {ENSTA Paris}, Institut Polytechnique de Paris, 91120 Palaiseau, France {\tt \footnotesize (andrea.simonetto@ensta-paris.fr)}. P. J. Goulart is with the Department of Engineering Science, University of Oxford, OX1 3PJ, United Kingdom ({\tt paul.goulart@eng.ox.ac.uk}).  
		This work was partially supported through the Government’s modern industrial strategy by Innovate UK under Project LEO (Ref. 104781).}}


\maketitle

\begin{abstract}
	We investigate both stationary and time-varying, nonmonotone generalized Nash equilibrium problems that exhibit symmetric interactions among the agents, which are known to be potential. As may happen in practical cases, however, we envision a scenario in which the formal expression of the underlying potential function is not available, and we design a semi-decentralized Nash equilibrium seeking algorithm. In the proposed two-layer scheme, a coordinator iteratively integrates possibly noisy and sporadic agents' feedback to learn the pseudo-gradients of the agents, and then design personalized incentives for them. On their side, the agents receive those personalized incentives, compute a solution to an extended game, and then return feedback measurements to the coordinator. In the stationary setting, our algorithm returns a Nash equilibrium in case the coordinator is endowed with standard learning policies, while it returns a Nash equilibrium up to a constant, yet adjustable, error in the time-varying case. As a motivating application, we consider the {ride-hailing} service provided by several competing companies with mobility as a service orchestration, necessary to both handle competition among firms and avoid traffic congestion.
\end{abstract}


\IEEEpeerreviewmaketitle

\section{Introduction}\label{sec:intro}
Noncooperative game theory represents a contemporary and pervasive paradigm for the modelling and optimization of modern multi-agent systems,
where agents are typically modelled as rational decision-makers that interact and selfishly compete for shared resources in a \emph{stationary} environment. Here, the (generalized) Nash equilibrium solution concept \cite{facchinei2007generalized} denotes a desired outcome of the game, which is typically self-learned by the agents through iterative procedures alternating distributed computation and communication steps \cite{salehisadaghiani2016distributed,salehisadaghiani2019distributed,ye2017distributed,gadjov2021exact}.	

Real-world scenarios, however, are rarely stationary. This fact, along with recent developments in machine learning and online optimization \cite{jadbabaie2015online,Bogunovic2016,shahrampour2017distributed,davis2019stochastic,dixit2019online,simonetto2019personalized,dall2020optimization}, has fostered the implementation of online multi-agent learning procedures, thus contributing in growing the interest for games where the population of agents ambitiously aim at tracking possibly \emph{time-varying} Nash equilibria online.  
In this paper we focus on both static and time-varying \emph{nonmonotone} \glspl{GNEP} that exhibit symmetric interactions among the agents, for which we design a semi-decentralized algorithm with convergence guarantees.
The presence of symmetric interactions, indeed, brings numerous advantages to the \gls{GNEP} that enjoys it. Among them, the underlying \gls{GNEP} is known to be potential with associated potential function \cite{ui2000shapley,la2016potential}, which implicitly entails the existence of a \gls{GNE}. Nonetheless, such a potential function frequently enables for the design of equilibrium seeking algorithms with convergence guarantees (especially in nonconvex setting \cite{heikkinen2006potential,fabiani2019multi,cenedese2019charging,fabiani2018distributed}).
However, unless one has a deep knowledge on the main quantities characterizing the symmetric interactions of the \gls{GNEP} at hand, finding the formal expression of the potential function is known to be a hard task \cite[Ch.~2]{la2016potential}\cite{hobbs2007nash}.
Moreover, in many in-network operations that require some degree of coordination it is highly desirable that the parameters of the agents' cost function, which reflect local sensitive data, stay private. 
In particular, our work is widely motivated by the {ride-hailing} application with \gls{maas} orchestration formally developed and described in \S \ref{sec:num_sim}.

In this framework, the proposed two-layer algorithm then reads as follows: in the outer loop, we endow a coordinator with an online learning procedure, aiming at iteratively integrating possibly noisy and sporadic agents' feedback to learn some of their private information, i.e., the pseudo-gradient mappings associated to the agents' cost functions. The reconstructed information is thereby exploited by the coordinator to design \emph{personalized incentives} \cite{simonetto2019personalized,ospina2020personalized,notarnicola2020distributed,Ospina2022} for the agents, which in turn compute a solution to an extended game, and then return (dis)satisfaction feedback measurements to the coordinator.

\subsection{Related work}
A recent research direction established the convergence of online distributed mirror descent-type algorithms in strictly monotone \glspl{GNEP} \cite{tampubolon2020coordinated}, aggregative games with estimated information \cite{tampubolon2019convergence}, or in price-based congestion control methods for generic noncooperative games \cite{tampubolon2020robust}.
Conversely, \cite{mertikopoulos2019learning} focused on the prediction of the long-term outcome of a monotone \gls{NEP},  also extended to the case of delays in the communication protocol \cite{zhou2018multi}. The convergence of no-regret learning policies with exponential weights in potential games within a (semi-)bandit framework was explored in \cite{cohen2017learning}, while \cite{cardoso2019competing} introduced an algorithm with sublinear Nash equilibrium regret under bandit feedback for time-varying matrix games, and \cite{duvocelle2018learning} showed that, in case of a slowly-varying monotone \gls{NEP}, the dynamic regret minimization allows the agents to track the sequence of~equilibria. Unlike the problem setting considered in this paper, \cite{fabiani2021nash} dealt with the quadratic class of stationary nonmonotone \glspl{GNEP} only where a feedback to the coordinator was provided at every iteration. 

We remark that our primary goal is to devise an algorithm to solve the \gls{GNE} seeking problem of an assigned noncooperative game. We do that by endowing a central coordinator with a learning procedure to reconstruct individual information of the agents, which is hence exploited to design personalized incentives, thus enabling for the equilibrium seeking. This is however misaligned with the overall goal of typical agents' utility and mechanism design, which mainly consists in designing games and payoff functions in which certain emergent behaviours coincide with a desirable outcome \cite{li2013designing,marden2018game,paccagnan2022utility}.

Finally, our semi-decentralized scheme may also be interpreted as a Stackelberg game in which the leader does not control any decision variable, albeit aims at minimizing the unknown potential function on the basis of optimistic conjectures on the followers' strategies \cite{kulkarni2015existence,fabiani2020local}. 
	
\subsection{Summary of contributions and paper organization}
In contrast to the aforementioned literature, we consider nonmonotone \glspl{GNEP} admitting symmetric interactions, for which we propose an online learning procedure based on \emph{personalized incentives} \cite{simonetto2019personalized,ospina2020personalized,notarnicola2020distributed}. The main contributions made can be summarized as follows:
\begin{itemize}
	\item We design a semi-decentralized scheme that allows the agents to compute (or track in a neighbourhood) a \gls{GNE} of a \emph{nonmonotone} \gls{GNEP} that admits an \emph{unknown} potential function (\S \ref{sec:prob_def}--\ref{sec:algorithm}). In particular:
	\begin{itemize}
		\item In the static case, we show that the proposed algorithm converges to a \gls{GNE} by exploiting the asymptotic consistency bounds characterizing typical learning procedures for the coordinator, such as \gls{ls} or \gls{gp} (\S \ref{sec:stat_case});
		\item In the time-varying setting, we show that the fixed point residual, our metric for assessing convergence, asymptotically behaves as $O(1)$, i.e., the proposed semi-decentralized scheme allows the agents to track a \gls{GNE} in a neighbourhood of adjustable size (\S \ref{sec:tv_case})
	\end{itemize}
	\item Inspired by real data available online \cite{nyc_ridehailing}, we develop a mathematical model capturing intrinsic features of the competition arising among companies of different size participating in the {ride-hailing} market with \gls{maas} orchestration. The developed model is then used as a case study to corroborate our theoretical results (\S \ref{sec:num_sim});
\end{itemize}
We also show that the design of the personalized incentives is key for the convergence of the algorithm, as they bring a twofold benefit: i) enabling the agents for the computation of a \gls{v-GNE} in the inner loop  by acting as a convexification terms for their cost functions; and ii) boosting the convergence and/or lessening the tracking error through a fine tuning of few parameters. 
Since our convergence results strongly rely on the knowledge of a private information held by each agent, specifically, the Lipschitz constant of their cost functions, in Appendix~\ref{sec:learning_lipschitz} we provide a possible learning-based solution to address this privacy issue, also accompanied by a dedicated analysis.
The proofs of theoretical results are all deferred to Appendix~\ref{sec:app_1}--\ref{sec:app_4}.

\vspace{-.5cm}
\subsection*{Notation}
$\N$, $\R$ and $\R_{\geq 0}$ denote the set of natural, real and nonnegative real numbers, respectively.  
$\bbS^{n}$ is the space of $n \times n$ symmetric matrices.
For vectors $v_1,\dots,v_N\in\mathbb{R}^n$ and $\mc I=\{1,\dots,N \}$, we denote $\bs{v} \coloneqq (v_1 ^\top,\dots ,v_N^\top )^\top = \mathrm{col}((v_i)_{i\in\mc I})$ and $ \bs{v}_{-i} \coloneqq \col(( v_j )_{j\in\mc I\setminus \{i\}})$. With a slight abuse of notation, we also use $\bs{v} = (v_i,\bs{v}_{-i})$.
$\mc{C}^1$ is the class of continuously differentiable functions. The mapping $F:\R^n \to \R^n$ is monotone on $\mc{X} \subseteq \R^n$ if $(F(x) - F(y))^\top(x - y) \, \geq  0$ for all $x, y \in \mc{X}$; strongly monotone if there exists a constant $c > 0$ such that $(F(x) - F(y))^\top(x - y) \geq c \|x - y\|^2$ for all $x, y \in \mc{X}$; hypomonotone if there exists a constant $c \geq 0$ such that $(F(x) - F(y))^\top(x - y) \geq -c \|x - y\|^2$ for all $x, y \in \mc{X}$. 
If $F$ is differentiable, $\mathrm{J}_F : \R^n \to \R^{n \times n}$ denotes its Jacobian matrix. 
Throughout the paper, variables with $t$ as subscript do not explicitly depend on time, as opposed when $t$ is an~argument.

\section{Problem formulation}\label{sec:prob_def}

We consider a noncooperative game $\Gamma \coloneqq (\mc{I}, (\mc{X}_i)_{i \in \mc{I}}, (g_i)_{i \in \mc{I}})$, with $N$ agents, indexed by the set $\mc{I} \coloneqq \{1, \ldots, N\}$. Each agent $i \in \mc{I}$ controls a local variable $x_i \in \mc{X}_i \subseteq \R^{n_i}$ and, at every discrete time instant $t \in \N$, aims at solving the following time-varying optimization problem:
\begin{equation}\label{eq:single_prob}
	\forall i \in \mc{I} : \left\{
	\begin{aligned}
		&\underset{x_i \in \mc{X}_i}{\textrm{min}} & & g_i (x_i,  \bs{x}_{-i}; t)\\
		&\hspace{.05cm}\textrm{ s.t. } & & h_i(x_i) + \sum_{j \in \mc{I} \setminus \{i\}} h_j(x_j)  \leq 0,
	\end{aligned}	
	\right.
\end{equation}
for some $g_i:\R^{n} \times \N \to \R$, $n = \sum_{i \in \mc{I}} n_i$, which denotes the private individual cost, whose value at time $t \in \N$ can be interpreted as the (dis)satisfaction of the $i$-th agent associated to the collective strategy $(x_i,  \bs{x}_{-i})$.
The collection of optimization problems in \eqref{eq:single_prob} amounts to a \gls{GNEP}, where every $h_i : \R^{n_i} \to \R^m$ is a map stacking $m$ coupling, yet locally separable, constraints among the agents. First, we define sets $\mc{X} \coloneqq \prod_{i \in \mc{I}} \mc{X}_i$ and $\mc{X}_i(\bs{x}_{-i}) \coloneqq \{x_i \in \mc{X}_i \mid h(x_i, \bs{x}_{-i}) \leq 0 \}$, with $h(x_i, \bs{x}_{-i}) \coloneqq h_i(x_i) + \sum_{j \in \mc{I} \setminus \{i\}} h_j(x_j)$, and then we introduce some standard assumptions:
\begin{standing}\label{standing:standard_assumptions}
	For each $i \in \mc{I}$, and for all $t \in \N$,
	\vspace{-.6cm}
	\begin{enumerate}
		\item[i)] The mapping $x_i \mapsto g_i(x_i, \bs{x}_{-i}; t)$ is of class $\mc{C}^1$ and has a $\ell_i$-Lipschitz continuous  gradient;
		\item[ii)] $\mc{X}_i$ is a nonempty, compact and convex set, $h_i : \R^{n_i} \to \R^m$ is a convex and of class $\mc{C}^1$ function.
		\QEDB
	\end{enumerate}
\end{standing}

The feasible set of the time-varying \gls{GNEP} $\Gamma$ thus coincides with $\Omega \coloneqq \{\bs{x} \in \mc{X} \mid h(\bs{x}) \leq 0\}$ \cite[\S 3.2]{facchinei2007generalized}. In the proposed time-varying context, we are then interested in designing an equilibrium seeking algorithm for the game $\Gamma$, according to the following popular definition of \gls{GNE}:

\begin{definition}\textup{(Generalized Nash equilibrium \cite{facchinei2007generalized})}\label{def:GNE}
	For all $t \in \N$, $\bs{x}^{\star}(t) \in \Omega$ is a \normaltext{\gls{GNE}} of the game $\Gamma$ if, for all $i \in \mc{I}$,
	\begin{equation}\label{eq:GNE}
		g_i(x^{\star}_i(t), \bs{x}^{\star}_{-i}(t);t) \leq \underset{y_i \in \mc{X}_i(\bs{x}^{\star}_{-i}(t))}{\normaltext{inf}} \; g_i(y_i, \bs{x}^{\star}_{-i}(t);t).
	\end{equation}
	\QEDB
\end{definition}

A collective vector of strategies $\bs{x}^{\star}(t) \in \Omega$ is therefore an equilibrium at time $t \in \N$ if no agent can decrease their objective function by changing unilaterally $\bs{x}^{\star}_{-i}(t)$ to any other feasible point. {Note that Definition~\ref{def:GNE} assumes a local connotation if one simply focuses on a certain neighbourhood of $\bs{x}^{\star}(t)$ for which condition \eqref{eq:GNE} holds true.}
Throughout the paper we will make use of the following assumption on the pseudo-gradient (or game) mapping $G: \R^n \times \N \to \R^n$, which is formally defined as $G(\bs{x};t) \coloneqq \col((\nabla_{x_i} g_i(x_i, \bs{x}_{-i};t))_{i \in \mc{I}})$:

\begin{standing}\label{standing:symmetry}
	For every $\bs{x} \in \Omega$ and $t \in \N$, $\mathrm{J}_G(\bs{x};t) \in \bbS^n$.
	\QEDB
\end{standing}

Roughly speaking, Standing Assumption~\ref{standing:symmetry} establishes that each pair of agents $(i,j) \in \mc{I}^2$ influences each other in an equivalent way. For the mapping $G$, this entails the existence of a differentiable, yet possibly \emph{unknown}, function $\theta : \R^n \times \N \to \R$ such that $G(\bs{x};t) = \nabla \theta(\bs{x};t)$, for all $\bs{x} \in \Omega$ and $t \in \N$ \cite[Th.~1.3.1]{facchinei2007finite}, which coincides with an {exact} potential function \cite{ui2000shapley,facchinei2011decomposition,la2016potential} for $\Gamma$ and can be characterized as stated next. 

\begin{lemma}\label{lemma:weak_conv}
	For all $t \in \N$, $\bs{x} \mapsto \nabla \theta(\bs{x};t)$ is $\ell$-Lipschitz continuous, while $\bs{x} \mapsto \theta(\bs{x};t)$ is $\ell$-weakly convex, i.e., $\bs{x} \mapsto \theta(\bs{x};t) + \tfrac{\ell}{2}\|x\|^2$ is convex, with $\ell \coloneqq \sum_{i \in \mc{I}}  \ell_i$.
	\QEDB
\end{lemma}

Note that $\theta$ is a smooth function that in principle may be nonconvex. Let $\Theta(t) \coloneqq \textrm{argmin}_{\bs{y} \in \Omega} \, \theta(\bs{y};t)$ be the set of its (local and global) constrained minimizers, assumed to be nonempty, and $\Theta^{\mathrm{s}}(t)$ be the set of its constrained stationary points, with $\Theta(t) \subseteq \Theta^{\mathrm{s}}(t)$, for all $t \in \N$. We stress that the nonemptiness of $\Theta(t)$ guarantees the existence of {an (at least local)} \gls{GNE} for $\Gamma$, since any $\bs{x}^\star(t) \in \Theta(t)$ satisfies the relation in \eqref{eq:GNE}. {It is well-known in potential game theory, indeed, that any minimum point of the exact potential function $\theta$ coincides with a \gls{GNE}, whose global or local nature depends on the point computed in $\Theta(t)$ -- see, for instance, \cite{roughgarden2010algorithmic,fabiani2019multi,cenedese2019charging}.}

\subsection{Main challenges and technical considerations}
In the considered framework, we identify three main critical issues
that rule out the possibility to compute a \gls{GNE} for the
\gls{GNEP} $\Gamma$ in \eqref{eq:single_prob} through standard arguments, thus fully motivating the design of a tailored learning procedure.

First, the time-varying nature of the optimization problems in \eqref{eq:single_prob} calls for an answer to the thorny question on whether there exist online learning policies that allow agents to track a Nash equilibrium over time (or to converge to one if the stage games stabilize). Even in the case of a potential game with known potential function, this is a challenging problem \cite{cohen2017learning}. 

In addition, despite the symmetry of interactions among agents, we note that for all $t \in \N$ the mapping $\bs{x} \mapsto G(\bs{x};t)$ may not be monotone, a key technical requirement for the most common solution algorithms for \glspl{GNEP} available in the literature, which compute a \gls{GNE} by relying on the (at least) monotonicity of the pseudo-gradient mapping \cite{salehisadaghiani2016distributed,salehisadaghiani2019distributed,ye2017distributed,gadjov2021exact}. 

Finally, we stress that Standing Assumption~\ref{standing:symmetry}, albeit quite mild and practically satisfied in several real-world scenarios \cite{rosenthal1973class,kattuman2004allocating,zhu2011multi,zhang2010cooperative}, is key to claim that the underlying \gls{GNEP} is potential.
However, unless one has a deep knowledge of the \gls{GNEP} at hand, finding the formal expression of the potential function is known to be a hard task \cite[Ch.~2]{la2016potential}. Thus, we assume not to have an expression for $\theta(\bs{x};t)$ that can be exploited directly for the equilibrium seeking algorithm design.

\begin{figure}[t!]
	\centering
	\ifTwoColumn
\else
	\includegraphics[width=.6\columnwidth]{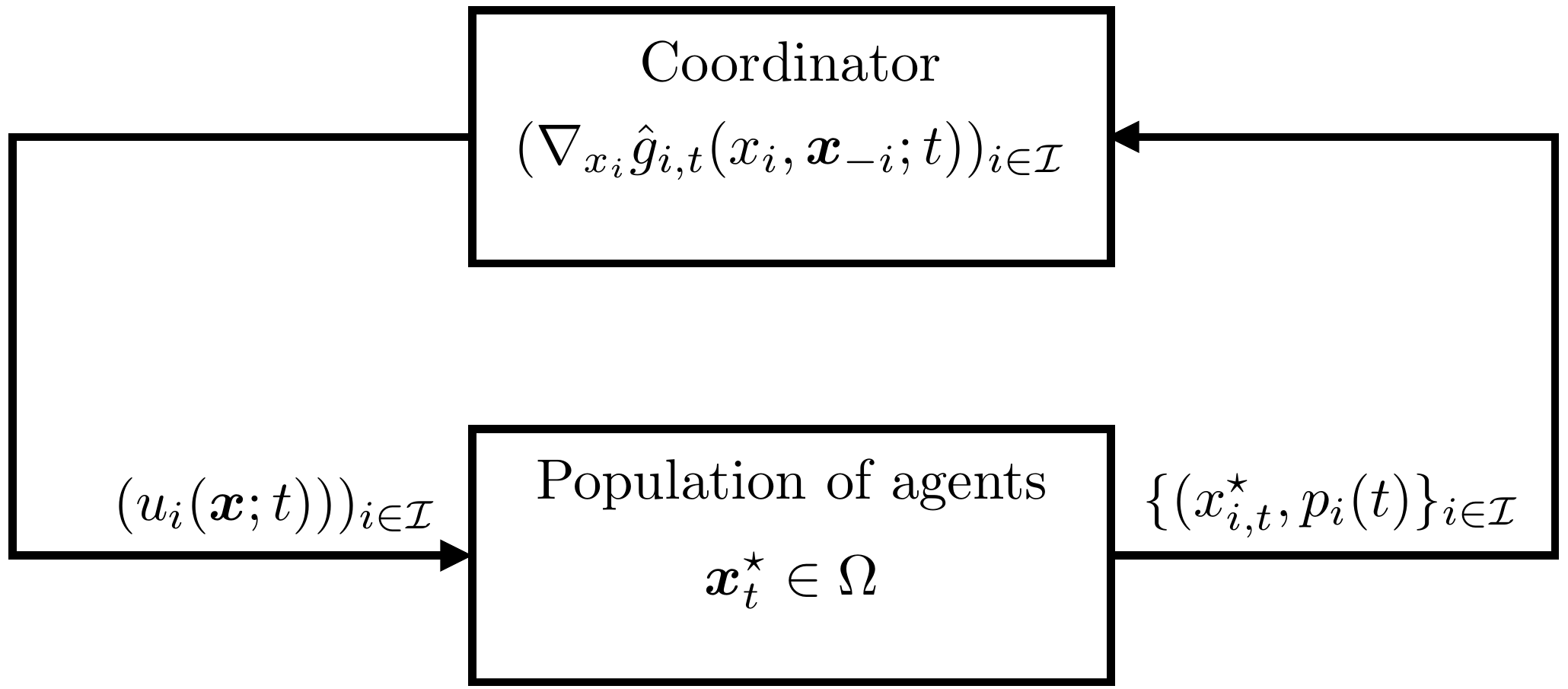}
	\fi	
	\caption{Personalized incentives as feedback design to steer the population to a point guaranteeing the ``minimum'' (dis)satisfaction, according to the unknown function $\theta$.}
	\label{fig:personalized_scheme}
\end{figure}

To address these crucial issues, we design \emph{personalized feedback functionals} $u_i : \R^{n} \times \N \to \R$ 	in the spirit of \cite{simonetto2019personalized,notarnicola2020distributed,ospina2020personalized,fabiani2021nash}, which are then used as ``control actions'' in the semi-decentralized scheme depicted in Fig.~\ref{fig:personalized_scheme}. Specifically, our goal is to steer the noncooperative agents to track minimizers of the unknown, time-varying function $\theta$, i.e., a \gls{GNE} of the game $\Gamma$, according to Definition~\ref{def:GNE}. Any $\bs{x}^\star(t) \in \Theta(t)$ can indeed be interpreted as a collective strategy that minimizes the (dis)satisfaction of the of agents, measured by the function~$\theta$.

\ifTwoColumn
\else
		\begin{algorithm}[!t]
	\caption{Two-layer semi-decentralized scheme}\label{alg:two_layer}
	\DontPrintSemicolon
	\SetArgSty{}
	\SetKwFor{ForAll}{for all}{do}{end forall}
	\smallskip
	\textbf{Iteration $(t \in \N)$:} \\
	\smallskip
	\begin{itemize}\setlength{\itemindent}{.2cm}
		\item[$\bullet$ (\texttt{S0})] Learn pseudo-gradients $(\nabla_{x_i} \hat{g}_{i, t - 1}(\bs{x}^\star_{t -1};t-1))_{i \in \mathcal{I}}$\\
		\smallskip
		\item[$\bullet$ (\texttt{S1})] Design personalized incentives $(u_i(\bs{x}; t))_{i \in \mc{I}}$\\
		\smallskip
		\item[$\circ$ (\texttt{S2})] Compute a \gls{GNE} of the extended game $\overline{\Gamma}$, $\bs{x}^\star_t \in \Omega$\\
		\smallskip
		\item[$\bullet$ (\texttt{S3})] Retrieve noisy agents' feedback $\{(x^\star_{i,t},p_{i}(t))\}_{i \in \mc{I}}$
	\end{itemize}
\end{algorithm}

\newpage

\section{Learning algorithm with\\ personalized incentives}\label{sec:algorithm}
\fi


\subsection{The two-layer algorithm}
The proposed approach is summarized in Algorithm~\ref{alg:two_layer}, where black-filled bullets refer to the tasks that have to be performed by a central coordinator, while the empty bullet to the one performed by the agents in $\mc{I}$. Thus, in the outer loop a central entity aims at learning online the unknown, time-varying function $\theta$ (or  its gradient mapping, $\nabla \theta$) by leveraging possibly noisy and sporadic agents' feedback on the private functions $g_i$'s (\texttt{S0}). On the basis of the estimated $\hat{g}_{i,t}$, at item (\texttt{S1}) the coordinator designs personalized incentive functionals $u_i${, which are successively communicated to the noncooperative agents taking part to the game. These latter then face with} an extended version of the \gls{GNEP} $\Gamma$ in \eqref{eq:single_prob} at item (\texttt{S2}), i.e., $\overline{\Gamma} \coloneqq (\mc{I}, (\mc{X}_i)_{i \in \mc{I}}, (f_i)_{i \in \mc{I}})$, with $g_i(x_i,  \bs{x}_{-i}; t) + u_i(x_i,  \bs{x}_{-i}; t) \eqqcolon f_i(x_i,  \bs{x}_{-i}; t)$ in place of $g_i(x_i,  \bs{x}_{-i}; t)$. Under a suitable choice of the personalized incentives, we will show that they act as regularization terms, as well as they trade-off convergence and robustness to the inexact knowledge of function $\theta$ and its gradient. Specifically, such incentives enable for the practical computation of an equilibrium of the extended game $\overline{\Gamma}$ at item (\texttt{S2}) through available solution algorithms for \glspl{GNEP} \cite{salehisadaghiani2016distributed,ye2017distributed,gadjov2021exact}, which typically require distributed computation and inner communication rounds among the noncooperative agents in $\mc I$.

Note that standard procedures in literature typically returns a \gls{v-GNE} \cite{cavazzuti2002nash,facchinei2007generalized}, which coincides to any solution to the \gls{GNEP} $\overline{\Gamma}$ that is also a solution to the associated \gls{VI}, i.e., any vector $\bs{x}_{t}^\star \in \Omega$ such that, for all $t \in \N$,	
\begin{equation}\label{eq:VI}
	(\bs{y} - \bs{x}^\star_t)^\top F(\bs{x}^\star_t; t) \geq 0, \, \text{ for all } \bs{y} \in \Omega, \, t \in \N,
\end{equation}	
where the mapping $F : \R^n \times \N \to \R^n$ is formally defined as $F(\bs{x}; t) \coloneqq \col((\nabla_{x_i} f_i(x_i, \bs{x}_{-i}; t))_{i \in \mc{I}}) = G(\bs{x};t) + U(\bs{x};t)$, and $U(\bs{x};t) \coloneqq \col((\nabla_{x_i} u_i(x_i, \bs{x}_{-i}; t))_{i \in \mc{I}})$. For these reasons, in referring to the computational step (\texttt{S2}), we tacitly assume that the agents compute a \gls{v-GNE} of the extended game $\overline{\Gamma}$.

Finally, at item (\texttt{S3}) the agents communicate feedback measures and their equilibrium strategies, $\{(x^\star_{i,t},p_{i}(t))\}_{i \in \mc{I}}$, with $p_{i}(t) \coloneqq g_i(\bs{x}^\star_t;t) + \varepsilon_{i, t}$, for some random variable $\varepsilon_i$, to the central entity, thus indicating to what extent the current equilibrium $\bs{x}^\star_t$ (dis)satisfies the entire population of agents.

Besides the time-varying nature of the original game $\Gamma$, introducing the personalized incentives $u_i(x_i,  \bs{x}_{-i}; t)$ further modifies $\Gamma$, thus forcing the agents to deal with an extended game $\overline \Gamma$ at each step (\text{S2}) of Algorithm~1. The relation between the equilibria of $\Gamma$ and $\overline \Gamma$ is hence not straightforward. What we prove later in the paper is that our approach allows the agents to compute (or track in a neighbourhood) a \gls{GNE} of the original \gls{GNEP} $\Gamma$, particularly that \gls{GNE} coinciding with a minimum of the unknown potential function $\theta(\bs x; t)$. We can therefore claim that, as $t \to \infty$, $\Gamma$ and $\overline \Gamma$ share at least one equilibrium, precisely that one coinciding with a minimum point of $\theta(\bs x; t)$. In different words, the learning procedure $\mathscr{L}$ characterizing the central coordinator, along with the parametric personalized incentives chosen, allows us to locally approximate the original nonmonotone GNEP $\Gamma$ around a certain GNE trajectory.

\subsection{Personalized incentives design}		
In view of Standing Assumption~\ref{standing:symmetry} we have
 $G(\bs{x}; t) = \nabla \theta(\bs{x}; t)$, and hence a natural approach to design the personalized incentives $u_i$ seems to iteratively learn and point a descent direction for the unknown function $\theta$, thus implicitly requiring one to estimate the pseudo-gradients $(\nabla_{x_i} g_{i}(x_i, \bs{x}_{-i};t))_{i \in \mc{I}}$, at every $t \in \N$.
Along the line of \cite{ospina2020personalized,fabiani2021nash}, we assume the central coordinator being endowed with a learning procedure $\mathscr{L}$ such that, at every outer iteration $t \in \N$ (Algorithm~\ref{alg:two_layer}, item (\texttt{S0})), it integrates the most recent agents' feedback $\{p_{i}(t-1)\}_{i \in \mc{I}}$ to  return an estimate of the pseudo-gradients, $(\nabla_{x_i} \hat{g}_{i, t-1}(\bs{x}^\star_{t-1}; t-1))_{i \in \mc{I}}$. 
A possible personalized incentive functional can hence be designed as	
\begin{equation}\label{eq:pers_feedback}
	u_i(\bs{x}; t) = \tfrac{1}{2} c({t}) \|x_i - x^+_{i}(t)\|^2, \text{ for all } i \in \mc{I},
\end{equation}	
where $x^+_{i}(t) \coloneqq x^\star_{i, t-1} + \xi(t) \nabla_{x_i} \hat{g}_{i, t-1}(\bs{x}^\star_{t-1}; t-1)$, for some parameters $c(t)$, $\xi(t)  \geq 0$, for all $t \in \N$. Unlike what one might expect, each $x^+_{i}(t)$ requires a positive sign for the gradient step $\xi(t) \nabla_{x_i} \hat{g}_{i, t-1}(\bs{x}^\star_{t-1}; t-1)$. However, note that this fact is not uncommon -- see, e.g., the recent Heavy Anchor method \cite[Eq.~(7)]{gadjov2021exact}. Moreover, we will also discuss later on how such a choice enables us to boost the convergence of Algorithm~\ref{alg:two_layer}
 or lessen the tracking error through a fine tuning of $\xi(t)$.

Thus, once the parametric form in \eqref{eq:pers_feedback} is fixed, we design suitable bounds for $c(t)$ and $\xi(t)$ in such a way that the sequence of \gls{GNE}, $(\bs{x}^\star_t)_{t \in \N}$, monotonically decreases $\theta(\bs{x};t)$ and converges to some point in $\Theta(t)$. As stressed in the previous section, the gain $c(t)$ is crucial to enable for the computation of a \gls{v-GNE} at item (\texttt{S2}) in Algorithm~\ref{alg:two_layer}, as formalized next:

\begin{proposition}\label{prop:strong_conv}
	Let $c(t) \geq 2 \ell$ for all $t \in \N$. Then, with the personalized incentives in \eqref{eq:pers_feedback}, the mapping $\bs{x} \mapsto F(\bs{x}; t)$ is $\ell$-strongly monotone, for all $t \in \N$.
	\QEDB
\end{proposition}

Thus, at every $t \in \N$, in \texttt{(S2)} the population of agents computes the (unique, see \cite[Th.~2.3.3]{facchinei2007finite}) \gls{v-GNE} associated to the extended version of the \gls{GNEP} in \eqref{eq:single_prob}, $\overline{\Gamma}$.
A key quantity for the convergence analysis of the proposed algorithm, both in the stationary and time-varying case, will be the fixed point residual $\Delta^\star_t \coloneqq \bs{x}^\star_{t} - \bs{x}^\star_{t-1}$, whose norm ``measures'' the distance to the points in $\Theta^{\mathrm{s}}$ when the function $\theta(\bs{x})$ is fixed in time.

\begin{lemma}\label{lemma:stationary}
	Let $(\bs{x}^\star_t)_{t \in \N}$ be the sequence of \normaltext{\gls{v-GNE}} generated by Algorithm~\ref{alg:two_layer} with $1-c(t)\xi(t)>0$, assume perfect reconstruction of the mapping $\bs{x} \mapsto G(\bs{x})$, and that $\|\Delta_t^\star\| = 0$ for some $\bs{x}^\star_t \in \Omega$. Then, $\bs{x}^\star_t \in \Theta^{\mathrm{s}}$.
	\QEDB
\end{lemma}	

As briefly seen,  choosing and tuning appropriately the parameters $c(t)$ and $\xi(t)$ is important to drive the sequence of \gls{v-GNE} $(\bs{x}^\star_t)_{t \in \N}$ along a descent direction for the unknown $\theta$, and ensuring $\|\Delta_t^\star\| \to 0$. Typically, we need $c(t) \geq 2\ell(t), 1-c(t)\xi(t)>0$, and we will see {in \S \ref{sec:num_sim}} how different choices could boost convergence and performance. In case of imperfect reconstruction of $\bs{x} \mapsto G(\bs{x})$, or in the time-varying setting, we also adopt the average value of $ \|\Delta_t^\star\|$ over a certain horizon of length $T \in \N$, i.e., $\frac{1}{T}\sum_{t\in \mc{T}} \|\Delta^\star_t\|$, $\mc{T}\coloneqq\{1, \ldots, T\}$, as a metric for the convergence of the sequence $(\bs{x}^\star_t)_{t \in \N}$ generated by Algorithm~\ref{alg:two_layer} to the stationary point set. 

We remark here that, on the one hand, finding the stationary points is the general goal in nonconvex setting \cite{scutari2017parallel}, and on the other hand, since Algorithm~\ref{alg:two_layer} generates monononically decreasing values for $\theta$, the application of simple perturbation techniques (e.g., \cite{escape}) can ensure that the stationary points to which we converge are in practice constrained local minima for $\theta$, namely points belonging to $\Theta \subseteq  \Theta^{\mathrm{s}}$, and therefore \gls{GNE} of the \gls{GNEP} $\Gamma$ in \eqref{eq:single_prob}, according to Definition~\ref{def:GNE}. {A technique that works well in practice is to use $G(\bs{x}) = \nabla \theta(\bs{x})$ to verify whether $\bs{x}$ is a minimum point of $\theta$ by small feasible perturbations, and if not, introduce that perturbation into \eqref{eq:pers_feedback}.}

\begin{remark}\label{remark:privacy_issue}
	The bounds on the parameters $c(t)$ and $\xi(t)$ provided in the paper assume the knowledge of the constant of weak convexity of $\theta$, $\ell = \sum_{i \in \mc{I}}  \ell_i$. As long as the coordinator is endowed with a learning policy, however, one may include this additional condition in the learning process, thus obtaining bounds that depend on $\hat{\ell} \geq 0$, the estimate of $\ell$. We discuss and elaborate more around this point in Appendix~\ref{sec:learning_lipschitz}.
	\QEDB
\end{remark}

\section{The stationary case}
\label{sec:stat_case}	
We start by discussing the case in which each $g_i$ in \eqref{eq:single_prob} is fixed in time, thus implying that $\theta(\bs{x};t) = \theta(\bs{x})$. First, we analyze the case of perfect reconstruction of the pseudo-gradient mappings $(\nabla_{x_i} \hat{g}_{i, {t-1}}(\bs{x}^\star_{t-1}))_{i \in \mc{I}}$ (\S \ref{sec:stat_case_1}), and then we investigate their inexact estimate (\S \ref{sec:stat_case_2}). 
Here, our result will be of the form  $(1/T)\sum_{t\in \mc{T}} \|\Delta^\star_t\| = O(1)$ in case the reconstruction error is non-vanishing. Otherwise, $(1/T)\sum_{t\in \mc{T}} \|\Delta^\star_t\| = 0$ (\S \ref{sec:stat_case_3}), thus recovering the results shown in \S \ref{sec:stat_case_1}.  	

\subsection{Online perfect reconstruction of the pseudo-gradients}\label{sec:stat_case_1}

In case the learning procedure $\mathscr{L}$ enables for $\nabla_{x_i} \hat{g}_{i, {t-1}}(\bs{x}^\star_{t-1}) = \nabla_{x_i} g_{i}(\bs{x}^\star_{t-1})$, $i \in \mc{I}$, by adopting the personalized incentives in \eqref{eq:pers_feedback} at every outer iteration $t \in \N$, we have the following result:

\begin{lemma}\label{lemma:desc_perfect}
	Let $c(t) \geq 2 \ell$ and $\xi(t)\in [0, 1/c(t))$, for all $t \in \N$. Then, with the personalized incentives in \eqref{eq:pers_feedback}, the vector $\Delta^\star_t$ is a descent direction for $\theta(\bs{x}^\star_{t-1})$, i.e., $\Delta^{\star^\top}_t \nabla \theta(\bs{x}^\star_{t-1}) < 0$.
	\QEDB
\end{lemma}

Then, if $c(t)$ (resp., $\xi(t)$) is large (small) enough, at every iteration $t \in \N$ of Algorithm~\ref{alg:two_layer}, the personalized functionals in \eqref{eq:pers_feedback} allow to point a descent direction for the unknown (dis)satisfaction function $\theta$. Next, we establish the convergence of the sequence of \gls{v-GNE} generated by Algorithm~\ref{alg:two_layer}.

\begin{proposition}\label{prop:convergence_perfect}
	Let $c(t) \geq 2 \ell$ and $\xi(t) \in [0, 1/c(t))$, for all $t \in \N$. With the personalized incentives in \eqref{eq:pers_feedback}, the sequence of \normaltext{\gls{v-GNE}} $(\bs{x}^\star_t)_{t \in \N}$, generated by Algorithm~\ref{alg:two_layer}, converges to some point in $\Theta^{\mathrm{s}}$.
	\QEDB
\end{proposition}

By introducing $\alpha(t) \coloneqq 1 - c(t) \xi(t)$, from the first step of the proof of Proposition~\ref{prop:convergence_perfect} we have that
$
\theta(\bs{x}^\star_t) \leq \theta(\bs{x}^\star_{t - 1}) - \ell (2 - \alpha(t))/2 \alpha(t) \|\Delta_t^\star\|^2
$, which points out that a fine tuning of the term $c(t)\xi(t)$ allows us to boost the convergence of Algorithm~\ref{alg:two_layer} to some point in $\Theta^{\mathrm{s}}$ (also observed on a numerical example in \cite[\S V]{fabiani2021nash}). This essentially explains the choice for a positive sign in the gradient step of \eqref{eq:pers_feedback}.
However, due to the presence of noise in the agents' feedback $\{p_i(t-1)\}_{i \in \mc{I}}$, it seems unlikely that the online algorithm $\mathscr{L}$ is able to return a perfect reconstruction of $(\nabla_{x_i} \hat{g}_{i, {t-1}}(\bs{x}^\star_{t-1}))_{i \in \mc{I}}$, at least at the beginning of the procedure in Algorithm~\ref{alg:two_layer}.

\subsection{Inexact estimate of the pseudo-gradients}\label{sec:stat_case_2}
At every outer iteration $t \in \N$, we assume the coordinator has available $K \in \N$ agents' feedback $\{p(k)\}_{k \in \mc{K}}$, $\mc{K} \coloneqq \{1, \ldots, K\}$, and $p(k) \coloneqq \col((p_{i}(k))_{i \in \mc{I}}) \in \R^{N}$, to estimate the gradients $(\nabla_{x_i} g_{i}(x_i, \bs{x}_{-i}))_{i \in \mc{I}}$ (and hence the mapping $\bs{x} \mapsto G(\bs{x})$).
The value of $K$ reflects situations in which the coordinator gathered information before starting the procedure ($K \geq t$), or it obtains sporadic feedback from the agents ($K < t$).
Without restriction, we make the following, standard assumption on the reconstructed mapping $\bs{x} \mapsto \hat{G}_{t}(\bs{x})$ directly, rather than on each single gradient \cite{dixit2019online,ospina2020personalized,dall2020optimization}.	

\begin{assumption}\label{ass:recons_error}
	For all $t \in \N$ and $\bs{x} \in \mc{X}$, $\hat{G}_t(\bs{x}) \coloneqq G(\bs{x}) + \epsilon_t$, and, for any $\delta_1 \in (0,1]$, there exists $\bar{t} < \infty$ and available $K_1$ agents' feedback such that 
	$$\prob\{\|\epsilon_{t}\| \leq \mathrm{e}(K) \mid \forall t \geq \bar{t}, \forall K\geq K_1\} \geq 1-\delta_1,$$ for some nonincreasing function $\mathrm{e} : \N \to \R_{\geq 0}$ such that $\mathrm{e}(K) < \infty$, for all $K \in \N$.
	\QEDB
\end{assumption}

With Assumption~\ref{ass:recons_error}, the reconstruction error on $\bs{x} \mapsto G(\bs{x})$ made by  $\mathscr{L}$ is bounded with high probability $1-\delta_1$ by some function of the available $K$ agents' feedback. 

\begin{remark}
Assumption~\ref{ass:recons_error} is reasonable for various learning strategies. Consider for instance a scalar \normaltext{\gls{ls}} estimator: take $e(K) = e$, i.e., a constant, then $\delta_1 \propto \mathrm{e}^{-\bar{t}}$, while for a function $e(K) \propto 1/t^{\vartheta}, \vartheta \in[0,1/2)$, then $\delta_1 \propto \Gamma(1/\varsigma, \bar{t}^\varsigma)$, with $\varsigma=1-2\vartheta\in(0,1]$, and $\Gamma(\cdot, \cdot)$ being the upper incomplete Gamma function. The latter is finite for our choice of parameters, and goes to zero as $\bar{t} \to \infty$. In both cases Assumption~\ref{ass:recons_error} is verified. See also~\normaltext{\cite[Lemma~A.4]{notarnicola2020distributed}} for the derivations and further extensions. \QEDB
\end{remark}

After defining quantities $\kappa(t) \coloneqq (1-\alpha(t))/2\alpha(t)$ and $\beta(t) \coloneqq \ell (2 - \alpha(t))/2 \alpha(t)$, we have the following result.	

\begin{lemma}\label{lemma:desc_perfect_recon}
	Let Assumption~\ref{ass:recons_error} hold true for some fixed $\delta_1 \in (0,1]$, $c(t) \geq 2 \ell$ and $\xi(t) \in [0, 1/c(t))$ for all $t \in \N$. Then, with the personalized incentives in \eqref{eq:pers_feedback}, for all $t \geq \bar{t}$ we have
	\begin{equation}\label{eq:descent_ball}
		\begin{aligned}
			\Delta^{\star^\top}_t \nabla \theta(\bs{x}^\star_{t-1}) &\leq \tfrac{\alpha(t) \kappa^2(t)}{\ell}  \mathrm{e}^2(K-1)- \left(\sqrt{\tfrac{\ell}{\alpha(t)}} \, \| \Delta_t^\star \| - \kappa(t) \sqrt{\tfrac{\alpha(t)}{\ell}} \, \mathrm{e}(K-1) \right)^2,
		\end{aligned}
	\end{equation} 
	with probability $1 - \delta_1$, for some $K \geq K_1$.
	\QEDB
\end{lemma}

In case of inexact estimate of the pseudo-gradients, the vector $\Delta_t^\star$ is not guaranteed to be a descent direction for the unknown function $\theta$.
In fact, the term $\mathrm{e}^2(K-1)$ rules out the possibility that the LHS in \eqref{eq:descent_ball} is strictly negative, albeit it can be made arbitrarily small through $\kappa(t)$ by an appropriate choice of the step-size $\xi(t)$. As in \S \ref{sec:stat_case_1}, the following bound characterizes the sequence of \gls{v-GNE} generated by Algorithm~\ref{alg:two_layer}.	

\begin{theorem}\label{th:conv_recons}
	Let Assumption~\ref{ass:recons_error} hold true for some fixed $\delta_1 \in (0,1]$, $c(t)\geq 2 \ell$ and $\xi(t) \in [0, 1/c(t))$, for all $t \in \N$. Moreover, let some $T \in \N$ be fixed, $\mc{T} \coloneqq \{\bar{t} + 1, \ldots, T + \bar{t}\}$ and, for any global minimizer $\bs{x}^\star \in \Theta$, $\Delta_{\bar{t}} \coloneqq \theta(\bs{x}^\star_{\bar{t}}) - \theta(\bs{x}^\star)$. Then, with the personalized incentives in \eqref{eq:pers_feedback}, the sequence of \normaltext{\gls{v-GNE}} $(\bs{x}^\star_t)_{t \in \mc{T}}$, generated by Algorithm~\ref{alg:two_layer}, satisfies the following relation with probability $1-\delta_1$
	\begin{equation}\label{eq:seq_recons}
		\begin{aligned}
			\tfrac{1}{T} \sum_{t\in \mc{T}} \|	\Delta_t^\star	\| &\leq \tfrac{1}{T \underline{\beta}} \sqrt{\sum_{t \in \mc{T}} \left( \beta(t) \Delta_{\bar{t}} + \tfrac{\bar{\beta} \kappa^2(t)}{\beta(t)} \mathrm{e}^2(q(t)) \right)}+ \tfrac{1}{T \underline{\beta}} \sum_{t \in \mc{T}}  \kappa(t) \mathrm{e}(q(t)).
		\end{aligned}
	\end{equation}
	Here, $\bar{\beta} \coloneqq \sum_{t \in \mc{T}} \beta(t)$, $\underline{\beta} \coloneqq \min_{t \in \mc{T}} \beta(t)$, and $q(t) \geq K_1$ is the number of available agents' feedback at the $t$-th outer iteration, $t \in \mc{T}$.
	\QEDB
\end{theorem}

Roughly speaking, Theorem~\ref{th:conv_recons} establishes that, with arbitrarily high probability, the average value of the residual $\|\Delta_t^\star\|$ over a certain horizon $T$ is bounded by the sum of two terms, which depend on the initial distance from a minimum for the unknown function $\theta$, and the reconstruction error $\mathrm{e}(\cdot)$. Note that the terms in the RHS can be made small by either choosing a small step-size $\xi(t)$, in order to make $\kappa(t)$ close to zero, or tuning the product $c(t) \xi(t)$ close to one, thus leading to a large $\underline{\beta}$. This latter choice, however, would increase the term involving the sub-optimal constant $\Delta_{\bar{t}}$, thus requiring an accurate trade-off in tuning the gain $c(t)$ and the step-size $\xi(t)$.
In the stationary case, to foster not exceedingly aggressive personalized actions the coordinator may then want to match the lower bound for $c(t)$, while striking a balance in choosing $\xi(t)$ to possibly boost the convergence of Algorithm~\ref{alg:two_layer}.

For simplicity, let us now assume that $\beta(t)$ is a constant term. From Assumption~\ref{ass:recons_error}, $\mathrm{e}(q(t)) \leq \mathrm{e}(K_1)$, and hence
$
\tfrac{1}{T} \sum_{t\in \mc{T}} \|	\Delta_t^\star	\| \leq O(1/\sqrt{T}) + O(1). 
$
We note that, as $T$ grows, $O(1/\sqrt{T})$ vanishes, and the average of $\|	\Delta_t^\star	\|$ stays in a ball whose radius depends on the number of agents' feedback $q(t)$ made available to perform (\texttt{S0}) in Algorithm~\ref{alg:two_layer} and, specifically, on the learning strategy $\mathscr{L}$. Next, we analyze the bound above under the lens of different learning procedures.

\subsection{Specifying the learning strategy $\mathscr{L}$}\label{sec:stat_case_3}
By requiring that the reconstruction error is bounded in probability, Assumption~\ref{ass:recons_error} is quite general and it holds true under standard assumptions for \gls{ls} and \gls{gp} approaches to learning $G$. In particular, we have the following:
\begin{itemize}
	\item In parametric learning, if $G(\bs{x})$ is modelled as an affine function of the learning parameters $\eta$'s, then setting up an \gls{ls} approach to minimize the loss between the model parameters $\eta$ and the agents' feedback leads to a convex quadratic program. Due to the large-scale properties of \gls{ls} (under standard assumptions), the error term $\mathrm{e}(q(t))$ behaves as a normal distribution, for which Assumption~\ref{ass:recons_error}  holds true (see~\cite[Lemma~A.4]{notarnicola2020distributed}), and $\lim_{K\to\infty} \mathrm{e}(K) = 0$.
	\item In non-parametric learning, suppose $G(\bs{x})$ is a sample path of a \gls{gp} with zero mean and a certain kernel. Due to the large-scale property of such regressor and under standard assumptions, also in this case Assumption~\ref{ass:recons_error}  holds true (see \cite{simonetto2019personalized}) and $\lim_{K\to\infty} \mathrm{e}(K) = 0$.
\end{itemize}	
Note that, in general, $q(t) \propto t$. Therefore, since $\sum_{t \in \mc{T}} \mathrm{e}(q(t)) = o(T)$, for the cases above we obtain
$
\lim_{T \to \infty} \tfrac{1}{T} \sum_{t\in \mc{T}} \|	\Delta_t^\star	\| = 0,
$
{
	thus implying that, for $T$ large enough, the average of $\|\Delta_t^\star	\|$ converges to $0$. If, in addition, there exists some $\bar K \ge 0$ such that $\mathrm{e}(K) = 0$ for all $K \ge \bar K$, then one is allowed to recover exactly the results obtained for the perfect reconstruction case shown in \S \ref{sec:stat_case_1}.
}

\section{The time-varying case}
\label{sec:tv_case}	
We now investigate the \gls{GNEP} in \eqref{eq:single_prob} in case the local cost function of each agent $g_i$ varies in time, thus implying that also the function $\theta$ is non-stationary. Our goal is still to design the parameters defining the personalized incentives to track a time-varying \gls{GNE} that minimizes the (dis)satisfaction function, i.e., some $\bs{x}^\star(t) \in \Theta(t)$, both in case of perfect (\S \ref{sec:tv_case_1}) and inexact reconstruction (\S \ref{sec:tv_case_2}) of the pseudo-gradient mapping. 

To start, we make the following typical assumptions in the literature on online optimization \cite{jadbabaie2015online,shahrampour2017distributed,dall2020optimization}.	
\begin{assumption}\label{ass:tv}
	For all $t \in \N$ and $\bs{x} \in \Omega$, it holds that
	\vspace{-.6cm}
	\begin{enumerate}
		\item[i)] $|\theta(\bs{x} ; t) - \theta(\bs{x} ; t-1) | \leq \mathrm{e}_{\theta}$, for $0 \leq \mathrm{e}_{\theta} < \infty$;
		\item[ii)] For all $i \in \mc{I}$, $\|\nabla_{x_i} g_{i}(\bs{x}; t) - \nabla_{x_i} g_{i}(\bs{x}; t-1)\| \leq \mathrm{e}_{\nabla_i}$, for $0 \leq \mathrm{e}_{\nabla_i} < \infty$;
		\item[iii)] $\|\bs{x}^{\star}(t) - \bs{x}^{\star}(t-1)\| \leq \mathrm{e}_{\delta}$, for $0 \leq \mathrm{e}_{\delta} < \infty$.
			\QEDB
	\end{enumerate}
\end{assumption}
Assumptions~\ref{ass:tv} i) and ii) essentially bound the variation in time of both the unknown function $\theta$ and the pseudo-gradient mappings, while Assumptions~\ref{ass:tv} iii) guarantees the boundedness of the distance between two consecutive minima such that $\bs{x}^{\star}(t) \in \Theta(t)$ and $\bs{x}^{\star}(t-1) \in \Theta(t-1)$. Note that, with these standard assumptions in place, an asymptotic error term of the form of $O(1)$ is inevitable~\cite{jadbabaie2015online,mokhtari2016online,NaLi2020,dall2020optimization}. 

\begin{lemma}\label{lemma:grad_tv}
	Let Assumption~\ref{ass:tv} ii) hold true. For all $t \in \N$, $\|\nabla \theta(\bs{x} ; t) - \nabla \theta(\bs{x} ; t-1) \| \leq \mathrm{e}_{\nabla}$, with $\mathrm{e}_{\nabla}  \coloneqq \sum_{i \in \mc{I}}  \mathrm{e}_{\nabla_i}$.
	\QEDB
\end{lemma}

\subsection{Online perfect reconstruction of the pseudo-gradients}
\label{sec:tv_case_1}	

In case the learning procedure $\mathscr{L}$ allows for $\nabla_{x_i} \hat{g}_{i, t - 1}(\bs{x}^\star_{t -1};t-1) = \nabla_{x_i} g_{i}(\bs{x}^\star_{t -1};t-1)$, for all $i \in \mc{I}$ and $t \in \N$, we have the following ancillary results:

\begin{lemma}\label{lemma:desc_perfect_tv}
	Let Assumption~\ref{ass:tv} ii) hold true, $c(t) \geq 2 \ell$ and $\xi(t) \in [0, 1/c(t))$, for all $t \in \N$. Then, with the personalized incentives in $\eqref{eq:pers_feedback}$, for all $t \in \N$ we have
	\begin{equation}\label{eq:descent_ball_tv}
		\begin{aligned}
			\Delta^{\star^\top}_t \nabla \theta(\bs{x}^\star_{t-1}; t-1) &\leq \tfrac{1}{4 \alpha(t) \ell}  \mathrm{e}^2_{\nabla}- \left(\sqrt{\tfrac{\ell}{\alpha(t)}} \, \| \Delta_t^\star \| - \tfrac{1}{2}\sqrt{\tfrac{1}{\alpha(t) \ell}} \, \mathrm{e}_{\nabla} \right)^2.
		\end{aligned}
	\end{equation} 
	\QEDB
\end{lemma}

As in the stationary case in \S \ref{sec:stat_case_2}, the vector $\Delta_t^\star$ is not guaranteed to be a descent direction for the unknown mapping $\bs{x} \mapsto \theta(\bs{x}; t-1)$ in the sense of Lemma~\ref{lemma:desc_perfect}. In fact, the error $\mathrm{e}_{\nabla}$, introduced because of the time-varying nature of the pseudo-gradients, excludes that the LHS in \eqref{eq:descent_ball} is strictly negative.
The following bound characterizes the sequence of \gls{v-GNE} originating from Algorithm~\ref{alg:two_layer} in case $\mathscr{L}$ allows for a perfect reconstruction of the time-varying mapping $\bs{x} \mapsto G(\bs{x}; t)$.

\begin{theorem}\label{th:conv_tv}
	Let Assumption~\ref{ass:tv} hold true, $c(t)\geq 2 \ell$ and $\xi(t) \in [0, 1/c(t))$, for all $t \in \N$. Moreover, let some $T \in \N$ be fixed, $\mc{T} \coloneqq \{1, \ldots, T \}$ and, for any global minimizer $\bs{x}^\star(0) \in \Theta(0)$, $\Delta_{0} \coloneqq |\theta(\bs{x}^\star_{0}; 0) - \theta(\bs{x}^\star(0); 0)|$. Then, with the personalized incentives in \eqref{eq:pers_feedback}, the sequence of \normaltext{\gls{v-GNE}} $(\bs{x}^\star_t)_{t \in \mc{T}}$, generated by Algorithm~\ref{alg:two_layer}, satisfies the following relation
	\begin{equation}\label{eq:seq_tv}
		\begin{aligned}
			&\tfrac{1}{T} \sum_{t \in \mc{T}} \|	\Delta_t^\star	\| \leq \tfrac{\mathrm{e}_{\nabla}}{T\underline{\beta}} \sum_{t \in \mc{T}} \tfrac{1}{2 \alpha(t)}+ \tfrac{1}{T \underline{\beta}} \sqrt{\sum_{t \in \mc{T}}   \left( \beta(t) \left(\Delta_0 + T \phi \right) + \tfrac{\bar{\beta}}{4 \alpha^2(t) \beta(t)} \mathrm{e}^2_\nabla \right)}.
		\end{aligned}
	\end{equation}
	with $\underline{\beta} \coloneqq \min_{t \in \mc{T}} \beta(t)$, $\phi \coloneqq 2 \mathrm{e}_{\theta} + \tfrac{\ell}{2} \mathrm{e}^2_{\delta}$ and $\bar{\beta} \coloneqq \sum_{t \in \mc{T}} \beta(t)$.~\QEDB
\end{theorem}

Theorem~\ref{th:conv_tv} says that the average of the residual $\|\Delta_t^\star\|$ over the horizon $T$ is bounded by the sum of two terms, which depend on the initial sub-optimality of a computed \gls{v-GNE} compared to a minimum for the unknown function $\theta$, and several bounds on the variations in time of $\theta$, $G$ and constrained minima postulated in Assumption~\ref{ass:tv} and Lemma~\ref{lemma:grad_tv}.  
In this case, the coordinator may reduce the error in the RHS by properly tuning the product $c(t) \xi(t)$ close to one, thus leading to a large $\underline{\beta}$, and hence possibly boosting the convergence of Algorithm~\ref{alg:two_layer}.
In fact, if the parameter $\beta(t)$ is fixed in time, we obtain
$
\tfrac{1}{T} \sum_{t\in \mc{T}} \|	\Delta_t^\star	\| \leq O(1).
$
This inequality ensures that $\|	\Delta_t^\star	\|$ will be always contained into a ball of constant radius, whose value can be adjusted through $c(t)$ and $\xi(t)$.

\subsection{Inexact estimate of the pseudo-gradients}\label{sec:tv_case_2}	
As in \S \ref{sec:stat_case_2}, we consider the case in which, due to possibly noisy agents' feedback, the learning procedure $\mathscr{L}$ does not allow a perfect reconstruction of each time-varying gradient $g_i$, $i \in \mc{I}$. 
First, we postulate the time-varying counterpart of Assumption~\ref{ass:recons_error}, and then we provide a preliminary result.	
\begin{assumption}\label{ass:recons_error_tv}
	For all $t \in \N$ and $\bs{x} \in \mc{X}$, $\hat{G}_t(\bs{x}; t) \coloneqq G(\bs{x}; t) + \epsilon_t$, and, for any $\delta_1 \in (0,1]$, there exists $\bar{t} < \infty$ and available $K_1$ agents' feedback such that 
	$$\prob\{\|\epsilon_{t}\| \leq \mathrm{e}(K) \mid \forall t \geq \bar{t}, \forall K\geq K_1\} \geq 1-\delta_1,$$ 
	for some nonincreasing function $\mathrm{e} : \N \to \R_{\geq 0}$ such that $\mathrm{e}(K) < \infty$, for all $K \in \N$.
	\QEDB
\end{assumption}	
\begin{lemma}\label{lemma:desc_imperfect_recon_tv}
	Let Assumption~\ref{ass:tv} and \ref{ass:recons_error_tv} hold true for some fixed $\delta_1 \in (0,1]$, $c(t) \geq 2 \ell$ and $\xi(t) \in [0,1/c(t))$, for all $t \in \N$. With the personalized incentives in $\eqref{eq:pers_feedback}$, for all $t \geq \bar{t}$ we have
	\begin{equation}\label{eq:descent_ball_imp_tv_}
		\begin{aligned}
			&\Delta^{\star^\top}_t  \nabla \theta(\bs{x}^\star_{t-1}; t -1) \leq  \tfrac{1}{4 \alpha(t) \ell} \sigma^2(K,t)- \left(\sqrt{\tfrac{\ell}{\alpha(t)}} \|\Delta_t^\star\| - \tfrac{1}{2\alpha(t)} \sqrt{\tfrac{\alpha(t)}{\ell}} \sigma(K,t) \right)^2,
		\end{aligned}
	\end{equation} 
	where $\sigma(K, t) \coloneqq \mathrm{e}_{\nabla} + (1-\alpha(t))\mathrm{e}(K-1)$, with probability $1 - \delta_1$, for some $K \geq K_1$.
	\QEDB
\end{lemma}
Along the same line drawn for the stationary case with inexact reconstruction, we now provide the following bound on the sequence of \gls{v-GNE}, $(\bs{x}^\star_t)_{t \in \mc{T}}$, generated by Algorithm~\ref{alg:two_layer}. Note the slight abuse of notation in defining $\Delta_{\bar{t}}$, which is different from the one in Theorem~\ref{th:conv_recons}.
\begin{theorem}\label{th:conv_recons_tv}
	Let Assumption~\ref{ass:tv} and \ref{ass:recons_error_tv} hold true for some fixed $\delta_1 \in (0,1]$, $c(t)\geq 2 \ell$ and $\xi(t) \in [0, 1/c(t))$, for all $t \in \N$. Moreover, let some $T \in \N$ be fixed, $\mc{T} \coloneqq \{\bar{t} + 1, \ldots, \bar{t} + T \}$ and, for any global minimizer $\bs{x}^\star(\bar{t}) \in \Theta(\bar{t})$, $\Delta_{\bar{t}} \coloneqq |\theta(\bs{x}^\star_{\bar{t}}; \bar{t}) - \theta(\bs{x}^\star(\bar{t}); \bar{t})|$. Then, with the personalized incentives in \eqref{eq:pers_feedback}, the sequence of \normaltext{\gls{v-GNE}} $(\bs{x}^\star_t)_{t \in \mc{T}}$, generated by Algorithm~\ref{alg:two_layer}, satisfies the following relation
	\begin{equation}\label{eq:seq_recons_tv}
		\begin{aligned}
			&\tfrac{1}{T} \sum_{t \in \mc{T}} \|	\Delta_t^\star\| \leq \tfrac{1}{T\underline{\beta}} \sum_{t \in \mc{T}} \tfrac{1}{2\alpha(t)}  \sigma(K,t)+ \tfrac{1}{T\underline{\beta}} \sqrt{\sum_{t \in \mc{T}} \left( \beta(t) \left(\Delta_{\bar{t}} + T \phi\right) + \tfrac{\bar{\beta}}{4 \alpha^2(t) \beta(t)} \sigma^2(K,t)\right)},
		\end{aligned}
	\end{equation}
	with probability $1-\delta_1$, where $\underline{\beta} \coloneqq \min_{t \in \mc{T}} \beta(t)$, $\bar{\beta} \coloneqq \sum_{t \in \mc{T}} \beta(t)$ and $q(t) \geq K_1$ is the number of agents' feedback at the $t$-th outer iteration.
	\QEDB
\end{theorem}
Also in this case, the average of the residual $\|\Delta_t^\star\|$ over the horizon $T$ is bounded by the sum of two terms, which depend, among the others, on the reconstruction error of the mapping $G$ and its variations in time. We note that the bound in the RHS of \eqref{eq:pers_feedback} can be adjusted through an accurate choice of the gain $c(t)$ and the step-size $\xi(t)$.
Specifically, choosing a small $\xi(t)$ reduces the reconstruction error, hidden in the variable $\sigma$, while setting $c(t) \xi(t)$ close to one induces a large value for $\underline{\beta}$ (and for $\beta(t)$ as well), thus possibly eliminating the second term under the square root of \eqref{eq:seq_recons_tv}, and the one outside.

For simplicity, let us now suppose that the parameters $c(t)$ and $\xi(t)$ of the personalized incentives in \eqref{eq:pers_feedback} are fixed in time, namely $\beta(t)$ is a constant term. From \eqref{eq:seq_recons_tv}, we note that
$
\tfrac{1}{T} \sum_{t\in \mc{T}} \|	\Delta_t^\star	\| \leq O(1). 
$
Due to the time-varying nature of the problem in question, also in this case the average residual $\|	\Delta_t^\star	\|$ can not vanish as $T$ grows, albeit the radius of the error ball can be reduced through a fine tuning of $c(t)$ and $\xi(t)$.

\subsection{Specifying the time-varying learning strategy $\mathscr{L}$}\label{subsec:tv_learning}
In a time-varying setting, one cannot expect $\mathrm{e}(K)$ to vanish in general, since the time variations in $\theta$ are not supposed to be asymptotically vanishing~\cite{jadbabaie2015online,dall2020optimization}. Popular learning approaches include \gls{ls} with forgetting factors \cite{Mateos09giannakis} and time-varying \gls{gp} \cite{Bogunovic2016}, for which we have $\lim_{K\to\infty} \mathrm{e}(K) = O(1)$.

\section{{Ride-hailing} with \gls{maas} orchestration}
\label{sec:num_sim}
With the growing business related to {ride-hailing}, a \gls{maas} coordination platform appears indispensable to contrast the traffic congestion due to the increasing number of vehicles dispatched on the road, while facilitating the competition among service providers \cite{pandey2019needs,diao2021impacts,fabiani2022stochastic}. In this section we develop a mathematical model capturing intrinsic features of the problem considered, and then use it to verify our theoretical findings.

\subsection{Mathematical model description}
We consider a scenario in which $N$ companies compete to put the most vehicles (a capped local resource,  $0 \leq x_i \leq \bar{x}_i$) on the road to attract the most customers. 
During the day, each company aims at maximizing its (time-varying) profit $P_i$, which is implicitly related to how many cars it could currently put on the road to meet customer needs, properly discounted to account for, e.g., refusals rates or the time of the day, $t$. To this end, bigger companies can be na\"ively induced to dispatch as many cars as they own. However, this may cause traffic congestion, thus reducing the quality of the service provided, and therefore lessen what the company can charge for each ride. 
In fact, by leveraging their own experience, those big companies may estimate how many cars actually get customers on top of the available $x_i$ as, e.g., a \emph{concave} function $a_{i}(t) x_i - b_{i}(t) x_i^2 \leq x_i$, with $a_{i}, b_{i} \geq 0$ tuned accordingly. Therefore, assuming the same fare $r_i(t)$ applies per average trip to each costumer, the profit function of the $i$-th lead company can read as
$
P_i(x_i; t) = r_i(t) (a_{i}(t) x_i - b_{i}(t) x_i^2).
$
On the other hand, the strategies of smaller companies are typically less affected by traffic congestions, since the quality service is generally worse in the sense that they can dispatch a little number of cars on the road. 
In this case, their direct experience may suggest that the number of cars that actually get customers on the available $x_i$ can be modelled as a \emph{convex} function $a_{i}(t) x_i + b_{i}(t) x_i^2 \leq x_i$, thus reflecting the fact that the larger the number of deployed cars, the larger the possibility to cover enough space to be attractive. The overall time-varying profit hence reads as
$
P_i(x_i; t) = r_i(t) (a_{i}(t) x_i + b_{i}(t) x_i^2).
$ 

In addition to the profit, however, the companies also incur in costs that have to be minimized and vary during the day, such as gas consumed or miles travelled, here represented by some $d_i(t) \ge 0$. By assuming, for instance, the same cost associated to each vehicle per average trip, we model the overall cost as
$C_i(x_i, \bs{x}_{-i}; t) = d_i(t) x_i + \sum_{j \in \mc{I} \setminus \{i\}} w_{i,j} (x_i - x_j)^2$,
where $(x_i - x_j)^2$ enables for competition among equally-sized companies, weighted with $w_{i,j} = w_{j,i} \geq 0$ to preserve symmetries.  

After splitting the hours of a day in intervals, enumerated in the set $\mc{N}$ according to the estimated average travel time of each costumer with no shared trips from \cite{nyc_ridehailing} (about $15$ minutes), at every $t \in \mc{N}$ each firm $i \in \mc{I}$ aims at solving the following mutually inter-dependent optimization problem:
\begin{equation}\label{eq:single_prob_maas}
	\forall i \in \mc{I} : \left\{
	\begin{aligned}
		&\underset{x_i \in [\underline{x}_i, \bar{x}_i]}{\textrm{min}} & &  C_i(x_i, \bs{x}_{-i}; t) - P_i(x_i;t)\\
		&\hspace{.3cm}\textrm{ s.t. } & & A \bs{x} \leq q,
	\end{aligned}	
	\right.
\end{equation}
where $A$ and $q$ collect constraints so that the total number of vehicles involved in the service is capped, i.e., $\bs{1}^\top \bs{x} \in [\underline{x}, \bar{x}]$, for some $0 \leq \underline{x} \leq \bar{x}$, and the fact that firm $i$ wants to improve its service compared to $j$, thus requiring $x_i \geq x_j$. 

In the proposed scenario where the firms exhibit symmetries in the mixed convex-concave cost functions\footnote{The adopted convex-concave costs can be thought as proxy for the actual objective that has to be minimized, which however can be complicated by considering, e.g., sum of exponential or logarithmic functions parametrized in the sensitive private quantities the coordinator aims at learning.}, the time-varying, unknown {exact} potential function for the \gls{GNEP} in \eqref{eq:single_prob_maas} is
$$
\theta(\bs{x}; t) \coloneqq \sum_{i \in \mc{I}} \left( - P_i(x_i;t)  + \sum\limits_{j \in \mc{I}, j < i} C_i(x_i, \bs{x}_{-i}; t) \right),
$$
and the \gls{maas} platform aims at coordinating the whole {ride-hailing} service while avoiding traffic congestion. This can be achieved, for instance, by imposing extra fees, incentives or restrictions to the companies, possibly according to their size and turnover. 
However, note that the parameters $a_{i}$ and $b_{i}$ affecting the cost function of each firm, which are hence key to drive its strategy, can not be disclosed to the \gls{maas} platform, since they represent sensitive information, as opposed to the incurred cost $C_i(x_i; t)$  that can be estimated directly, as it depends on mileage, fuel consumption and number of deployed cars. {These data are indeed not sensitive, and therefore can be either publicly available as in \cite{nyc_ridehailing} (e.g., ride-hailing trips per day, average vehicles per day, daily trips per vehicle, minutes per trip, fare-box per trip, etc.), or may be estimated from them.}
Thus, 
a possible strategy requires the \gls{maas} platform to learn those time-varying parameters by leveraging feedback collected from users, e.g., on the price they are charged $r_i(t)$, and then design tailored incentives for the coordination. The fact that only a subset of deployed vehicles gets customers intrinsically represents the noise in the costumers' feedback, as the profit of each firm does not na\"ively coincide with the applied fare $r_i(t)$ times the number of vehicles $x_i(t)$.

\subsection{Numerical simulations}

\begin{table}
	\caption{Simulation parameters}
	\label{tab:parameters_perfect_rec}
	\centering\scalebox{0.95}{
	\begin{tabular}{lll}
		\toprule
		Parameter  & Description   & Value \\
		\midrule
		$\ell(t)$ & Constant of weak convexity & $[0.38, 5.2]$\\
		$a_i(t)$ & Discount parameter -- profit function & $[0.92, 0.94]$\\
		$b_i(t)$ & Discount parameter -- profit function & $[-0.41, 6.6] \!\times\! 10^{-4}$\\
		$r_i(t)$ & Trip fare & $[14, 30]$ [\$]\\
		$d_i(t)$ & Trip cost & $[4, 10]$ [\$]\\
		$P$ & Weight for the extragradient algorithm & $I$ \\
		$l_F$ & Horizontal scale -- \gls{gp} method  & $2 \times 10^3$ \\
		$\sigma_F$ & Vertical scale -- \gls{gp} method  & $10^4$ \\
		$\gamma$ & Discount factor -- \gls{gp} method & $0.8$ \\
		\bottomrule
	\end{tabular}}
\end{table}

\begin{figure}[t!]
	\centering
	\ifTwoColumn
\else
	\includegraphics[width=.6\columnwidth]{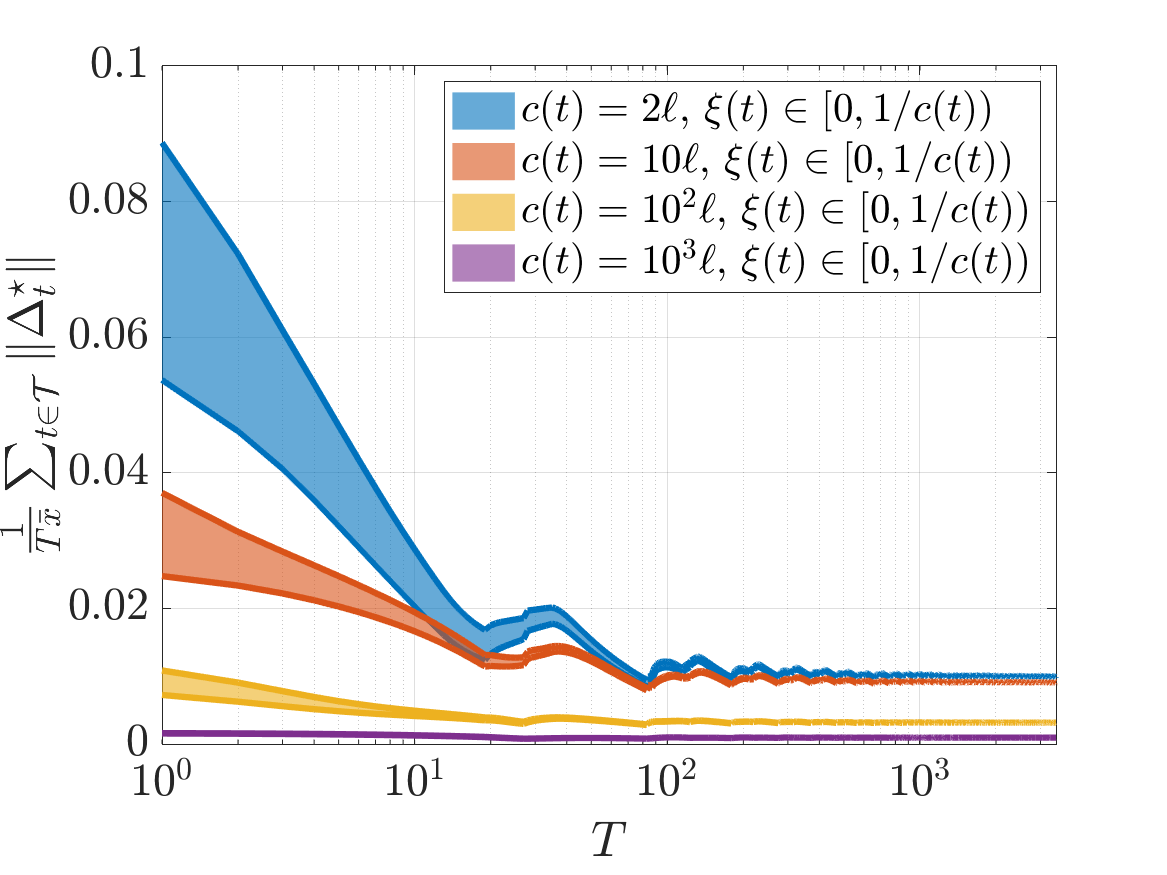}
	\fi	
	\caption{Residual $\Delta^\star_t$ over the time horizon $T$, normalized by the total number of deployable cars at any $t$, for different values of personalized incentive feedback gains $c(t)$ and $\xi(t)$ in \eqref{eq:pers_feedback}.}
	\label{fig:path_length_GP}
\end{figure}

The open data collected in New York City in April 2019 \cite{nyc_ridehailing} provide us information on $N = 5$ main companies: Yellow taxi, Uber, Lyft, Juno and Via. We stress that the quality of the {ride-hailing} service, measured as the total number of vehicles deployed on the road, coincides with an integer variable, i.e., $n_i = 1$, for all $i \in \mc{I} = \{1, \ldots, 5\}$, thus leading to a mixed-integer setting. However, since the fleet dimension of each firm we are considering is in the order of few thousands of vehicles (i.e., Juno and Via), or tens of thousands for bigger companies (Uber, Lyft, Yellow taxi), we consider a relaxed version by treating $x_i$ as a scalar continuous variable, and then rounding its value \cite{pandey2019needs}. For this reason, we roughly estimate a round-off error in the order of $10^{-3}$ for any \gls{GNE} computed at item \texttt{(S2)} in Algorithm~\ref{alg:two_layer} through an extragradient type method \cite{solodov1996modified} (thus neglecting the multi-agent nature of the inner loop). 
	
To learn the unknown terms characterizing each pseudo-gradient associated to the cost function of each company, i.e., the time-varying parameters $a_i$ and $b_i$, we assume the \gls{maas} platform being endowed with a \gls{gp} learning algorithm. Further numerical examples with other regression strategies can be found in the preliminary work \cite{fabiani2021nash}, where we considered a specific class of stationary nonmonotone \glspl{GNEP} only.

\begin{figure}[t!]
	\centering
	\ifTwoColumn
\else
	\includegraphics[width=.6\columnwidth]{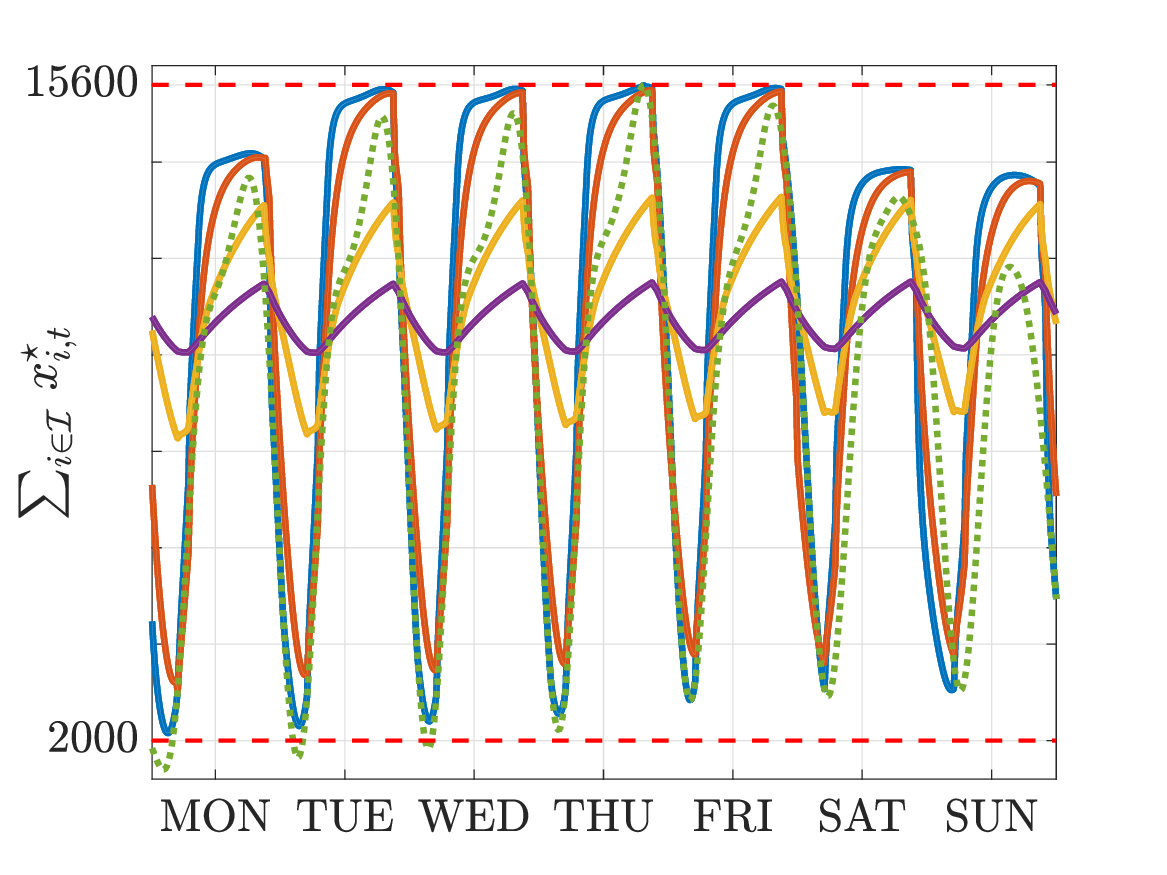}
	\fi	
	\caption{Overall service provided by the five firms over a whole week. The colour line follows the legend in Fig.~\ref{fig:path_length_GP}, while the dotted green line corresponds to the costumer requests obtained from data after a spline interpolation. The dashed red lines coincide with $\underline{x}$ and $\bar{x}$, respectively.}
	\label{fig:constraints_satisf}
\end{figure}

	The main parameters adopted to run numerical experiments are summarized in Table~\ref{tab:parameters_perfect_rec}, where $d_i(t)$ and $r_i(t)$ are time-varying functions designed to capture the costumer requests variability over an entire week (see the dotted green line in Fig.~\ref{fig:constraints_satisf} for instance). We stress that the information on the costumer demands is not directly exploited anywhere, except for the design of functions $d_i(t)$ and $r_i(t)$, i.e., the numerical results are obtained without artificially imposing the service requests as a reference trajectory.
	According to the data available at \cite{nyc_ridehailing}, in April 2019 Yellow taxi represented approximately the 27\% of the whole car-based mobility market, whereas among the {ride-hailing} firms Uber impacted for the 72\% on the market, thus representing the leading company with potentially no competitors, Lyft for the 19.7\%, Via and Juno for the 4.8\% and 3.4\%, respectively.
	Thus, by comparing the vehicle fleet size characterizing each firm, it turns out that the weight matrix $W \coloneqq [w_{i,j}]_{i \in \mc{I}, j \in \mc{I}}$ links Yellow taxi and Lyft, and Via and Juno only. These companies, indeed, compete each other to provide at least a comparable service for the sake of reputation.

\begin{figure}[t!]
	\centering
	\ifTwoColumn
\else
	\includegraphics[width=.6\columnwidth]{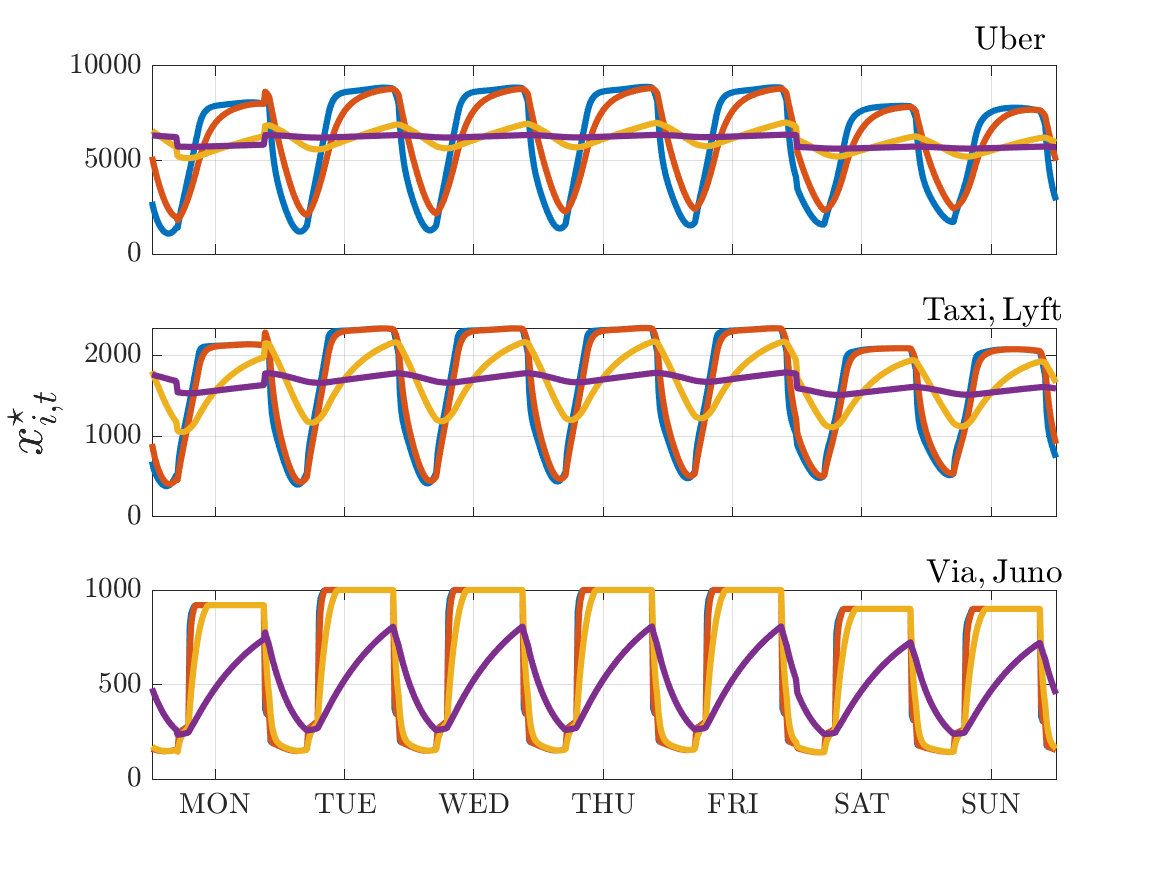}
	\fi	
	\caption{Optimal strategies pursued by the five companies over a week, coinciding with the number of vehicles deployed at every time interval. The colour line follows the legend in Fig.~\ref{fig:path_length_GP}.}
	\label{fig:weekly_strategy}
\end{figure}

\begin{figure}[t!]
	\centering
	\ifTwoColumn
\else
	\includegraphics[width=.6\columnwidth]{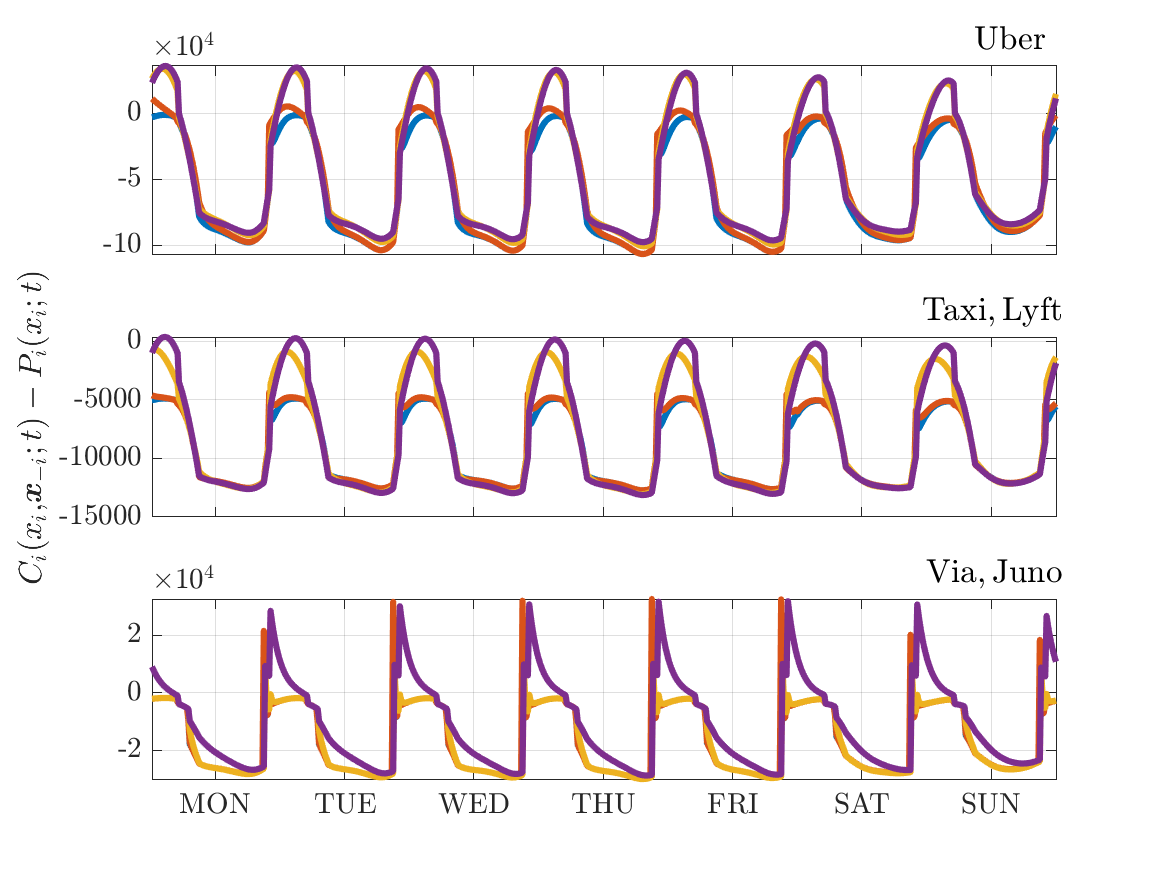}
	\fi	
	\caption{Profit/cost experienced by the companies over a week. The colour line follows the legend in Fig.~\ref{fig:path_length_GP}.}
	\label{fig:weekly_profit}
\end{figure}

	In Fig.~\ref{fig:path_length_GP} is reported the average value of the residual $\Delta^\star_t$ over the time horizon $T$, normalized by its magnitude. This latter, indeed, was identified as a candidate metric to assess convergence of the proposed algorithm, with different values of personalized incentive gains $c(t)$ and $\xi(t)$.  Note that increasing $\xi(t)$ noticeably reduces the residual $\Delta^\star_t$ for short time horizon $T$, while a larger $c(t)$ seems providing a smaller asymptotic error. Given the time variability of the costumer demands, however, this may result in firms' optimal strategies leading to an overall unsatisfactorily service, as evidenced by Fig.~\ref{fig:weekly_strategy} and, more prominently, by Fig.~\ref{fig:constraints_satisf}. In fact, for $c(t) = 2 \ell(t)$ (or even $10 \ell(t)$, blue and brown lines), the \gls{maas} platform manages to accomplish the task of serving hundreds of thousands of costumers per week, whereas for larger values of $c(t)$ the companies do not capture the variability of the demand (yellow and violet lines), adopting almost constant strategies over the week and hence experiencing higher costs, as reported in Fig.~\ref{fig:weekly_profit}. This represents the trade-off that one has to strike between performance of Algorithm~\ref{alg:two_layer} and expected/desired behaviour of the companies and, more importantly, of the overall service provided. For this reason, the \gls{maas} platform is a key tool enabling for competition among companies in {ride-hailing} mobility while, at the same time, guaranteeing a certain degree of service by satisfying costumer requests.

\begin{figure}[t!]
	\centering
	\ifTwoColumn
\else
	\includegraphics[width=.6\columnwidth]{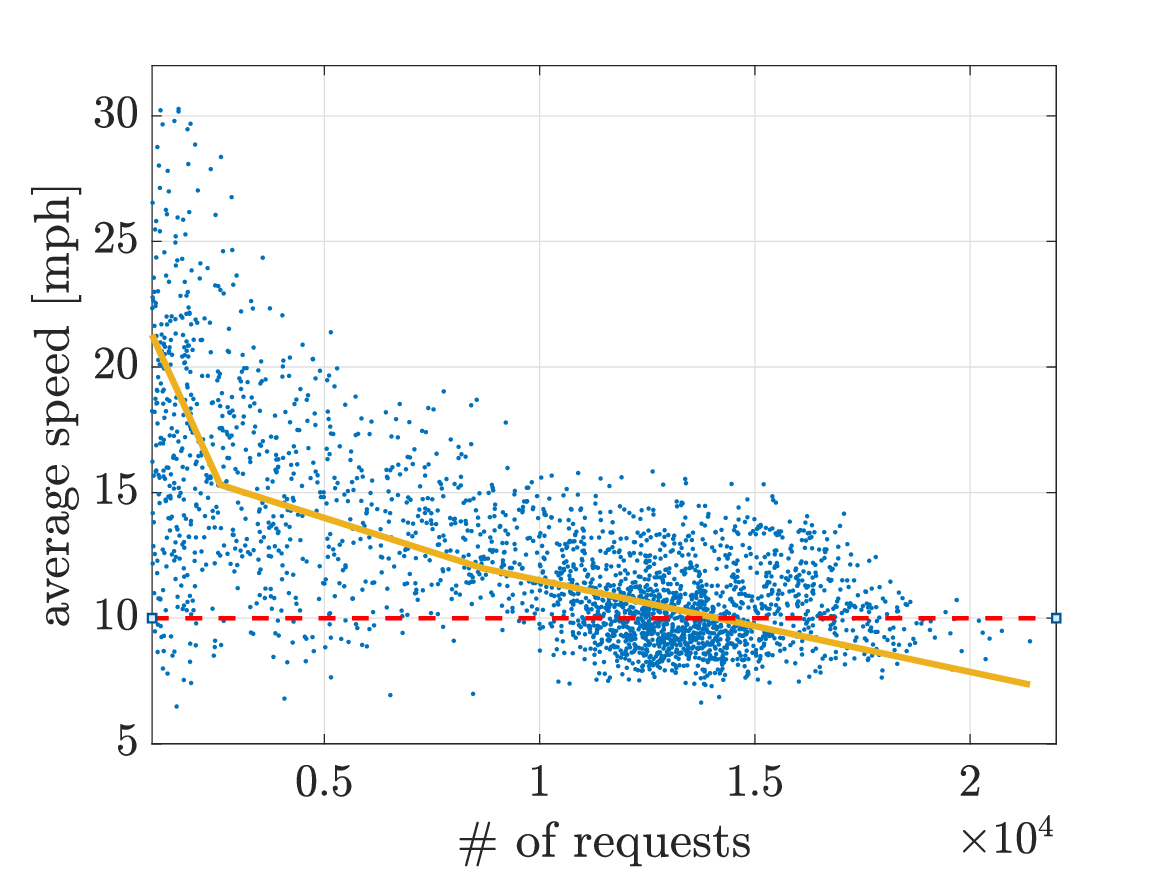}
	\fi
	\caption{Piecewise-affine interpolation linking the number of costumer requests and average speed per trip.}
	\label{fig:interp_requests_speed}
\end{figure}

\begin{figure}[t!]
	\centering
	\ifTwoColumn
\else
	\includegraphics[width=.6\columnwidth]{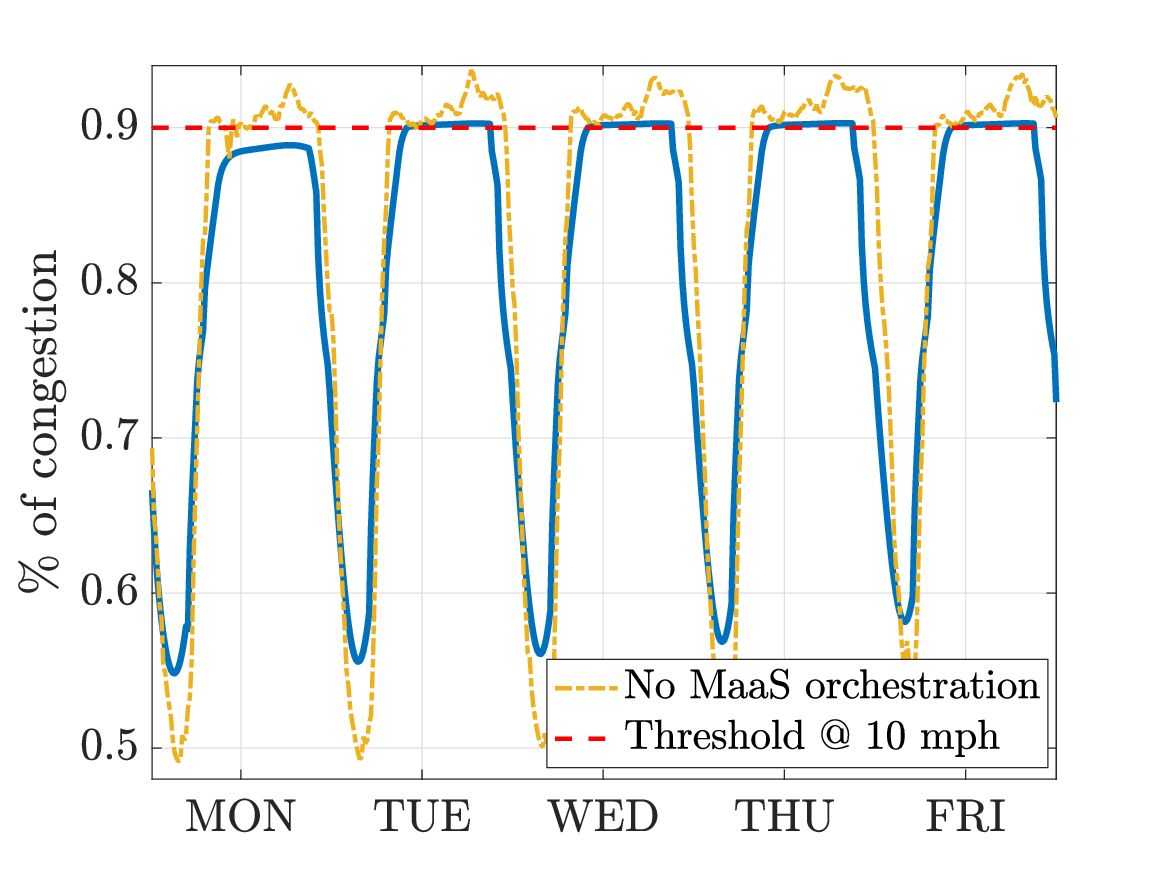}
	\fi
	\caption{Comparison between the traffic congestion caused by ridehailing service with and without \gls{maas} orchestration.}
	\label{fig:weekly_congestion}
\end{figure}

We finally give further motivation to the importance of a \gls{maas} orchestration system as the one proposed in this paper in terms of potential impact to alleviate traffic congestion. In \cite{nyc_ridehailing}, indeed, is shown that each vehicle contributes to two rides per hour on average, while for the rest of the time is idle. By considering an average trip length of 15 minutes, this potentially leads to more than 50\% of cars deployed on the road than necessary (i.e., in accordance with the costumer requests), thus heavily affecting the urban traffic, especially in big cities such as New York. Then, by leveraging available data, we first extrapolate the relation among number of costumer requests and average speed per trip, as described in Fig.~\ref{fig:interp_requests_speed} where we have adopted a piece-wise affine function to characterize this relation. Successively, by assuming a conservative estimate of up to the 25\% of vehicles dispatched more than the actual number of costumer requests, we evaluate the effect of the \gls{maas} platform on the traffic congestion in Fig.~\ref{fig:weekly_congestion}. By focusing on weekdays, in which the traffic congestion produces the greatest costumer discomfort, the \gls{maas} orchestration seems to represent a crucial component to avoid deploying an excessive number of cars on the road, thus avoiding to reach and surpass the threshold at 90\% corresponding to the critical speed per trip, identified as 10 miles per hour.

\section{Conclusion and Outlook}
\label{sec:conclusion}

We have shown that a suitable design of parametric personalized incentives is crucial to compute (or track in a neighbourhood) \gls{GNE} in nonmonotone \gls{GNEP} characterized by symmetric interactions among agents, both in static and time-varying setting. First, the designed functionals act as regularization terms of the agents' cost functions, thus allowing the agents for the practical computation of a \gls{v-GNE} at each outer iteration of the proposed two-layer algorithm. Then, they provide a mean to boost the convergence of the algorithm in the stationary case, or to adjust the asymptotic error obtained in a time-varying setting. In the static case, the proposed algorithm converges to a \gls{GNE} by exploiting the asymptotic consistency bounds characterizing standard learning procedures for the coordinator, such as \gls{ls} or \gls{gp} , while in a time-varying setting our metric for assessing convergence asymptotically behaves as $O(1)$, and hence the semi-decentralized scheme allows the agents to track a \gls{GNE} in a neighbourhood of adjustable size.

Future research directions may include for instance
the design of other possible parametric personalized incentives,
as well as non-parametric ones, in which the pseudo-gradient mappings are directly learned through, e.g., neural networks. The case of possibly approximate computation of a \gls{v-GNE} at item \texttt{(S2)} may be explored, also working towards the relaxation of the symmetric interaction requirement.

\appendix

\subsection{Proofs of \S \ref{sec:prob_def}}\label{sec:app_1}
\textit{Proof of Lemma~\ref{lemma:weak_conv}}: 	Let some $t \in \N$ be fixed. Then, in view of Standing Assumption~\ref{standing:standard_assumptions} and \ref{standing:symmetry}, for any $\bs{x}$, $\bs{y} \in \R^n$, the following chain of inequalities hold
$
	\|	\nabla \theta(\bs{x};t) - \nabla \theta(\bs{y};t)	\| = \|	 G(\bs{x};t) -  G(\bs{y};t)	\| = \|	\col((\nabla_{x_i} g_i(x_i, \bs{x}_{-i};t) - \nabla_{y_i} g_i(y_i, \bs{y}_{-i};t)_{i \in \mc{I}}) \| \leq \textstyle\sum_{i \in \mc{I}} \| \nabla_{x_i} g_i(x_i, \bs{x}_{-i};t) - \nabla_{y_i} g_i(y_i, \bs{y}_{-i};t) \| \leq \textstyle\sum_{i \in \mc{I}} \ell_i \| \bs{x} - \bs{y} \| \ \eqqcolon \ell \| \bs{x} - \bs{y} \|.
$
Finally, the fact that $\bs{x} \mapsto \theta(\bs{x}; t)$ is a $\mc{C}^1$-smooth function with $\ell$-Lipschitz continuous gradient directly entails that, for any $\bs{x} \in \Omega$, the auxiliary function $\psi(\bs{x};t) \coloneqq \theta(\bs{x};t) + \tfrac{\ell}{2} \| \bs{x} \|^2$ is convex, which in turn implies the $\ell$-weak convexity of $\bs{x} \mapsto \theta(\bs{x};t)$ \cite{davis2019stochastic,mai2020convergence}.
\hfill\qedsymbol

\subsection{Proofs of \S \ref{sec:algorithm}}\label{sec:app_2}
{
	\textit{Proof of Proposition~\ref{prop:strong_conv}}: With the personalized incentive functionals in \eqref{eq:pers_feedback},  $F(\bs{x}; t) = G(\bs{x}; t) + c(t) (\bs{x} - \bs{x}^\star_{t-1} - \xi(t)\hat{G}_{t-1}(\bs{x}^\star_{t-1}))$. Thus, for any $\bs{x}$, $\bs{y} \in \Omega$, and $t \in \N$, we have that:
	$
		(\bs{x} - \bs{y})^\top (F(\bs{x};t) - F(\bs{y};t)) = (\bs{x} - \bs{y})^\top (G(\bs{x}; t) + c(t) (\bs{x} - \bs{x}^\star_{t-1} - \xi(t)\hat{G}_{t-1}(\bs{x}^\star_{t-1})) - G(\bs{y}; t) - c(t) (\bs{y} - \bs{x}^\star_{t-1} - \xi(t)\hat{G}_{t-1}(\bs{x}^\star_{t-1}))) = (\bs{x} - \bs{y})^\top (G(\bs{x}; t) + c(t) \bs{x} - G(\bs{y}; t) - c(t) \bs{y}) = (\bs{x} - \bs{y})^\top (\nabla\theta(\bs{x}; t) + c(t) \bs{x} - \nabla\theta(\bs{y}; t) - c(t) \bs{y}),
	$
	where this last equality follows from Standing Assumption~\ref{standing:symmetry}. Then, let us introduce the auxiliary function $\psi(\bs{x};t) \coloneqq \theta(\bs{x}; t) + \tfrac{c(t)}{2} \| \bs{x} \|^2$. Note that, in case $c(t) \geq 2 \ell$ for all $t \in \N$, $\bs{x} \mapsto \psi(\bs{x};t)$ is $\ell$-strongly convex in view of the $\ell$-weak convexity of $\theta$ proved in Lemma~\ref{lemma:weak_conv}, which yields to
	$$
	\begin{aligned}
		&(\bs{x} - \bs{y})^\top (F(\bs{x};t) - F(\bs{y};t)) \\  
		&\qquad=(\bs{x} - \bs{y})^\top (\nabla\theta(\bs{x}; t) + c(t) \bs{x} - \nabla\theta(\bs{y}; t) - c(t) \bs{y})\\
		&\qquad= (\bs{x} - \bs{y})^\top (\nabla\psi(\bs{x}; t) - \nabla\psi(\bs{y}; t)) \geq \ell \|\bs{x} - \bs{y}\|^2,
	\end{aligned}
	$$
	i.e., the definition of strongly monotone mapping.
	\hfill\qedsymbol
}

\smallskip

\textit{Proof of Lemma~\ref{lemma:stationary}}: Consider the stationary points of $\theta(\bs{x})$ as the ones satisfying the first-order optimality conditions: $w \nabla\theta(\bs{x}) + \mathcal{N}_{\Omega}(\bs{x}) \ni \bs{0}$, where $w>0$ is any scalar scaling of the cost $\theta$, and $ \mathcal{N}_{\Omega}$ is the normal cone operator of the feasible set $\Omega$. At every $t \in \N$, $\bs{x}^\star_t$ with perfect reconstruction is the solution to 
$
\nabla\theta(\bs{x}^\star_t) + c(t) (\bs{x}^\star_t - \bs{x}^\star_{t-1} - \xi(t) \nabla \theta(\bs{x}^\star_{t-1})) + \mathcal{N}_{\Omega}(\bs{x}^\star_t) \ni \bs{0},
$
If $\|\Delta^\star_t\| = 0$ then $\bs{x}^\star_t = \bs{x}^\star_{t-1}$, which means that $\bs{x}^\star_t$ is the solution to 
$
(1- c(t)\xi(t)) \nabla \theta(\bs{x}^\star_{t}) + \mathcal{N}_{\Omega}(\bs{x}^\star_t) \ni \bs{0},
$
which satisfies the first-order optimality conditions for $\theta(\bs{x})$ and it is therefore one of its stationary points. 
\hfill\qedsymbol

\subsection{Proofs of \S \ref{sec:stat_case}}\label{sec:app_3}

\textit{Proof of Lemma~\ref{lemma:desc_perfect}}: 	In view of item (\texttt{S2}), at every outer iteration $t \in \N$ we have by \eqref{eq:VI} that $\bs{x}^\star_t$  satisfies $(\bs{y} - \bs{x}^\star_t)^\top F(\bs{x}^\star_t; t) \geq 0$ for all $\bs{y} \in \Omega$, and therefore, since $\bs{x}^\star_{t-1} \in \Omega$ as it is a \gls{v-GNE} at $t - 1$, $\Delta_t^{\star^\top} F(\bs{x}^\star_t; t) \leq 0$. Thus, by adding and subtracting the term $\Delta_t^{\star^\top} G(\bs{x}^\star_{t-1})$, we obtain
$
\Delta_t^{\star^\top} G(\bs{x}^\star_{t-1}) \leq \Delta_t^{\star^\top} (G(\bs{x}^\star_{t-1}) - (G(\bs{x}^\star_t) + U(\bs{x}^\star_t; t))).
$
With the personalized functionals in \eqref{eq:pers_feedback} and the perfect estimate of the pseudo-gradients, we have that $U(\bs{x}^\star_t; t) = c(t) (\Delta_t^\star - \xi(t) G(\bs{x}^\star_{t-1}))$. Then, in view of Standing Assumption~\ref{standing:symmetry}, it follows that 
$$
\begin{aligned}
	(1 - c(t)\xi(t)) \Delta_t^{\star^\top} \nabla \theta(\bs{x}^\star_{t-1}) &\leq \Delta_t^{\star^\top} ( \nabla \theta(\bs{x}^\star_{t-1}) + c(t) \bs{x}^\star_{t-1} - \nabla \theta(\bs{x}^\star_t) - c(t) \bs{x}^\star_{t}).
\end{aligned}
$$
Now, by defining $\alpha(t) \coloneqq 1 - c(t) \xi(t)$ for all $t \in \N$, let us introduce the auxiliary function $\psi(\bs{x};t) \coloneqq \theta(\bs{x}) + \tfrac{c(t)}{2} \| \bs{x} \|^2$. Note that, in case $c(t) \geq 2 \ell$ for all $t \in \N$, $\bs{x} \mapsto \psi(\bs{x};t)$ is $\ell$-strongly convex in view of the $\ell$-weak convexity of $\theta$, proved in Lemma~\ref{lemma:weak_conv}. As such, $\Delta_t^{\star^\top} ( -\nabla \psi(\bs{x}^\star_{t-1}; t) + \nabla \psi(\bs{x}^\star_t;t)) \geq \ell \| \Delta_t^\star \|^2$. Therefore, we have
$
	\alpha(t) \Delta_t^{\star^\top} \nabla \theta(\bs{x}^\star_{t-1}) \leq \Delta_t^{\star^\top} ( \nabla \psi(\bs{x}^\star_{t-1}; t) - \nabla \psi(\bs{x}^\star_t;t)) \leq -\ell \| \Delta_t^\star \|^2.
$
Finally, we obtain the desired result, namely $\Delta^{\star^\top}_t \nabla \theta(\bs{x}^\star_{t-1}) \leq - (\ell/\alpha(t)) \| \Delta_t^\star \|^2 < 0$, by imposing $\alpha(t) = 1 - c(t) \xi(t) > 0$, i.e., $\xi_{i} < 1/c(t) \leq 1/2\ell$, for all $t \in \N$.
\hfill\qedsymbol

\smallskip

\textit{Proof of Proposition~\ref{prop:convergence_perfect}}: 	By combining the descent lemma \cite[Prop.~A.24]{bertsekas1997nonlinear} and Lemma~\ref{lemma:desc_perfect}, the sequence $(\bs{x}^\star_t)_{t \in \N}$ satisfies
$$
	\begin{aligned}
		\theta(\bs{x}^\star_t) &\leq \theta(\bs{x}^\star_{t - 1}) + (\Delta_t^\star)^\top \nabla \theta(\bs{x}^\star_{t - 1}) + \tfrac{\ell}{2} \|\Delta_t^\star\|^2\\
		&\leq \theta(\bs{x}^\star_{t - 1}) - \ell \tfrac{2 - \alpha(t)}{2 \alpha(t)} \|\Delta_t^\star\|^2.
	\end{aligned}
$$
Thus, by imposing $(2 - \alpha(t))/2 \alpha(t) > 0$, which entails that $0 \leq \xi(t) < 1/c(t) \leq 1/2\ell$, then the sequence $(\theta(\bs{x}^\star_t))_{t \in \N}$ shall converge to a finite value, as $\theta(\bs{x}^\star_t) \to -\infty$ can not happen in view of the compactness of $\mc{X}$ (Standing Assumption~\ref{standing:standard_assumptions}). Therefore, by the continuity of $\theta$, the convergence of $(\theta(\bs{x}^\star_t))_{t \in \N}$ entails that $\mathrm{lim}_{t \to \infty} \, \|\Delta_t^\star\| = 0$, and hence by Lemma~\ref{lemma:stationary}, the bounded sequence of (feasible by definition) points $(\bs{x}^\star_t)_{t \in \N} \to \bs{x} \in \Theta^{\mathrm{s}}$.
\hfill\qedsymbol

\smallskip

\textit{Proof of Lemma~\ref{lemma:desc_perfect_recon}}:	By mimicking the same steps at the beginning of the proof of Lemma~\ref{lemma:desc_perfect}, in view of Assumption~\ref{ass:recons_error} we obtain, for all $t \geq \bar{t}$,
$
		\alpha(t) \Delta_t^{\star^\top} \nabla \theta(\bs{x}^\star_{t-1}) \leq \Delta_t^{\star^\top}  ( \nabla \theta(\bs{x}^\star_{t-1}) + c(t) \bs{x}^\star_{t-1} - \nabla \theta(\bs{x}^\star_t) - c(t) \bs{x}^\star_{t}) + c(t) \xi(t) \Delta_t^{\star^\top}  \epsilon_{t-1}.
$
In case $c(t) \geq 2 \ell > 0$ for all $t \in \N$, the first term in the RHS is always upper bound by $ -\ell \| \Delta_t^\star \|^2$ in view of the $\ell$-strongly convexity of the auxiliary function $\theta(\bs{x}) + \tfrac{c(t)}{2} \| \bs{x} \|^2$. On the other hand, the second term attains its maximum positive module when the two vectors $\Delta_t^\star$ and $\epsilon_{t-1}$ are aligned (since $\xi(t) \geq 0$), thus leading to the following chain of inequalities
$$
\begin{aligned}
	\alpha(t) \Delta_t^{\star^\top} \nabla \theta(\bs{x}^\star_{t-1}) &\leq  -\ell \| \Delta_t^\star \|^2 + c(t) \xi(t) \|\Delta_t^\star\|\| \epsilon_{t-1}\|\\
	&\leq -\ell \| \Delta_t^\star \|^2 + c(t) \xi(t) \|\Delta_t^\star\| \mathrm{e}(K-1),
\end{aligned}
$$
where the last relation follows from Assumption~\ref{ass:recons_error}, for any $K \geq K_1$, with probability $1 - \delta_1$. Thus, for $\alpha(t) > 0$, i.e., $\xi(t) < 1/c(t)$, and $c(t) \xi(t) = 1 - \alpha(t)$, we obtain
\begin{equation}\label{eq:descent}
	\Delta_t^{\star^\top} \nabla \theta(\bs{x}^\star_{t-1}) \leq  - \tfrac{\ell}{\alpha(t)} \| \Delta_t^\star \|^2 + \tfrac{1-\alpha(t)}{\alpha(t)} \| \Delta_t^\star \| \mathrm{e}(K-1).
\end{equation}
Finally, adding and subtracting the term $((1-\alpha(t))^2/4\alpha(t)\ell) \mathrm{e}^2(K-1)$ in the RHS directly yields~to
$$
\begin{aligned}
	\Delta_t^{\star^\top} \nabla \theta(\bs{x}^\star_{t-1}) \leq &-\left(\sqrt{\tfrac{\ell}{\alpha(t)}} \, \| \Delta_t^\star \| - \tfrac{1-\alpha(t)}{2\alpha(t)}\sqrt{\tfrac{\alpha(t)}{\ell}} \, \mathrm{e}(K-1) \right)^2 + \tfrac{(1-\alpha(t))^2}{4\alpha(t)\ell}  \mathrm{e}^2(K-1),
\end{aligned}
$$
with probability $1-\delta_1$. The proof concludes by substituting $\kappa(t) \coloneqq (1-\alpha(t))/2 \alpha(t)$ in the inequality above.
\hfill\qedsymbol

\smallskip

\textit{Proof of Theorem~\ref{th:conv_recons}}: As in the proof of Proposition~\ref{prop:convergence_perfect}, we start by combining the descent lemma \cite[Prop.~A.24]{bertsekas1997nonlinear} and Lemma~\ref{lemma:desc_perfect_recon} (specifically, the relation in \eqref{eq:descent}). In view of Assumption~\ref{ass:recons_error}, for all $t \geq \bar{t}$ and $K \geq K_1$, we obtain
$$
\theta(\bs{x}^\star_t) \leq \theta(\bs{x}^\star_{t - 1})   - \ell \tfrac{2 - \alpha(t)}{2 \alpha(t)} \|\Delta_t^\star\|^2 + \tfrac{1 - \alpha(t)}{\alpha(t)} \|\Delta_t^\star\| \mathrm{e}(K-1)
$$	
with probability $1 - \delta_1$. By defining $\beta(t) \coloneqq \ell (2 - \alpha(t))/2 \alpha(t)$, which is guaranteed to be strictly greater than zero in case $\xi(t) \in [0, 1/c(t))$, for all $t \in \N$ (and hence for $t \geq \bar{t}$), and focusing on the second and third terms in the RHS of this latter inequality, we can complete the square by adding and subtracting $((1-\alpha(t))^2/4\alpha^2(t) \beta(t)) \mathrm{e}^2(K-1)$. After few algebraic manipulations, we obtain the following inequality
\begin{equation}\label{eq:sample_complexity_structure}
	\begin{aligned}
		\theta(\bs{x}^\star_t) &\leq \theta(\bs{x}^\star_{t - 1}) -\beta(t) \left(	\|\Delta_t^\star\| - \tfrac{1-\alpha(t)}{2 \alpha(t) \beta(t)} \mathrm{e}(K-1) \right)^2+ \tfrac{(1-\alpha(t))^2}{4 \alpha^2(t) \beta(t)} \mathrm{e}^2(K-1).
	\end{aligned}
\end{equation}

From now on, we follow the proofline introduced in \cite[Th.~3.1]{davis2019stochastic}, \cite[Th.~1]{mai2020convergence}, as in \eqref{eq:sample_complexity_structure} we have obtained a typical structure characterizing the sample complexity analysis for sample-based stochastic algorithms (see, for instance, \cite[Eq.(1.3)]{davis2019stochastic}). Specifically, by fixing some $T \in \N$ and summing up the inequality above over $t \in \mc{T} \coloneqq \{\bar{t} + 1, \ldots, T + \bar{t}\}$, in view of the possibly sporadic communication between the agents and the central coordinator, let $q(t)$ be a scalar indicating the number of agents' feedback adopted by the learning procedure $\mathscr{L}$ at the $t$-th outer iteration. In particular, by assuming to initialize Algorithm~\ref{alg:two_layer} with $q(0) = K_1$ agents' feedback, $q(t) = q(t-1) + 1$ if, at the $t$-th iteration, the agents communicate $p(t)$ to the coordinator, $q(t) = q(t-1)$ otherwise. Then, for any global minimizer $\bs{x}^\star \in \Theta$, let $\kappa(t) \coloneqq (1-\alpha(t))/2 \alpha(t)$. With probability $1-\delta_1$ we obtain
$
	\theta(\bs{x}^\star) \leq \theta(\bs{x}^\star_T) \leq \theta(\bs{x}^\star_{\bar{t}}) - \sum_{t \in \mc{T}} \beta(t) (	\|\Delta_t^\star\| -  \tfrac{\kappa(t)}{\beta(t)} \mathrm{e}(q(t)) )^2 + \sum_{t \in \mc{T}} \tfrac{\kappa^2(t)}{\beta(t)} \mathrm{e}^2(q(t)).
$
After moving the term $\theta(\bs{x}^\star)$ to the RHS, the term with $\beta(t)$ to the LHS, and introducing $\Delta_{\bar{t}} \coloneqq \theta(\bs{x}^\star_{\bar{t}}) - \theta(\bs{x}^\star) \geq 0$, then:
\begin{equation}\label{eq:1}
	\sum_{t \in \mc{T}} \beta(t) \left(	\|\Delta_t^\star\| \!-\! \tfrac{\kappa(t)}{\beta(t)} \mathrm{e}(q(t)) \right)^2 \!\leq\! \Delta_{\bar{t}}  \!+\! \sum_{t \in \mc{T}} \tfrac{\kappa^2(t)}{\beta(t)} \mathrm{e}^2(q(t)).
\end{equation}

Since $\beta(t) >0$ for all $t \in \mc{T}$, we can leverage the Jensen's inequality \cite[Th.~3.4]{rockafellar1970convex} on the convex function $(\cdot)^2$ to lower bound the summation in the LHS. 
To this end, we normalize the coefficient by multiplying and dividing by $\bar{\beta} \coloneqq \sum_{t \in \mc{T}} \beta(t)$, and we define $\hat{\beta}(t) \coloneqq \beta(t)/ \bar{\beta}$, thus obtaining 
$
	\bar{\beta} (\sum_{t \in \mc{T}} \hat{\beta}(t)  (	\|\Delta_t^\star\| - \tfrac{\kappa(t)}{\beta(t)} \mathrm{e}(q(t)) )^2 ) \geq \bar{\beta} (\sum_{t \in \mc{T}} \hat{\beta}(t)	(\|\Delta_t^\star\| - \tfrac{\kappa(t)}{\beta(t)} \mathrm{e}(q(t)) ) )^2.
$
Then, replacing this latter inequality into \eqref{eq:1}, 
$$
\begin{aligned}
	\bar{\beta} \left(\sum_{t \in \mc{T}} \hat{\beta}(t)	\left(\|\Delta_t^\star\| - \tfrac{\kappa(t)}{\beta(t)} \mathrm{e}(q(t)) \right) \right)^2 &= \tfrac{1}{\bar{\beta}} \left(\sum_{t \in \mc{T}} {\beta}(t)	\left(\|\Delta_t^\star\| - \tfrac{\kappa(t)}{\beta(t)} \mathrm{e}(q(t)) \right) \right)^2\\
	&\leq \Delta_{\bar{t}}  + \sum_{t \in \mc{T}} \tfrac{\kappa^2(t)}{\beta(t)} \mathrm{e}^2(q(t)),
\end{aligned}
$$
and performing few algebraic manipulations directly yields to,	
$$
\begin{aligned}
	&\sum_{t \in \mc{T}} \beta(t) \left(\|\Delta_t^\star\| - \tfrac{\kappa(t)}{\beta(t)} \mathrm{e}(q(t)) \right) \leq\sqrt{\sum_{t \in \mc{T}} \left(\beta(t) \Delta_{\bar{t}} + \tfrac{\bar{\beta} \kappa^2(t)}{\beta(t)} \mathrm{e}^2(q(t)) \right)}.	
\end{aligned}
$$

By bringing the term $- \sum_{t \in \mc{T}} \kappa(t) \mathrm{e}(q(t))$ in the RHS, note that the obtained inequality still holds true if we premultiply both sides by $1/T$, thus obtaining the average over $T$.
Moreover, we have that $(1/T) \sum_{t \in \mc{T}} \beta(t) \|\Delta_t^\star\| \geq (1/T) \underline{\beta} \sum_{t \in \mc{T}} \|\Delta_t^\star\|$, with $\underline{\beta} \coloneqq \textrm{min}_{t \in \mc{T}} \beta(t)$, which directly yields to the relation in \eqref{eq:seq_recons} with probability $1-\delta_1$, hence concluding the proof.
\hfill\qedsymbol

\subsection{Proofs of \S \ref{sec:tv_case}}\label{sec:app_4}

\textit{Proof of Lemma~\ref{lemma:grad_tv}}: 	Let some $t \in \N$ be fixed. Then, in view of Standing Assumption~\ref{standing:standard_assumptions}, \ref{standing:symmetry}, and Assumption~\ref{ass:tv} ii), for any $\bs{x} \in \R^n$ the following chain of inequalities hold:
$
	\|	\nabla \theta(\bs{x};t) - \nabla \theta(\bs{x};t-1)	\| = \|	 G(\bs{x};t) -  G(\bs{x};t-1)	\| \leq \sum_{i \in \mc{I}} \| \nabla_{x_i} g_i(\bs{x};t) - \nabla_{x_i} g_i(\bs{x};t-1) \| \leq \sum_{i \in \mc{I}} \mathrm{e}_{\nabla_i} \eqqcolon \mathrm{e}_{\nabla}.
$
\hfill\qedsymbol

\smallskip

\textit{Proof of Lemma~\ref{lemma:desc_perfect_tv}}: 	In view of item (\texttt{S2}) in Algorithm~\ref{alg:two_layer}, at every outer iteration $t \in \N$ we have by equation \eqref{eq:VI} that $\bs{x}^\star_t$  satisfies $(\bs{y} - \bs{x}^\star_t)^\top F(\bs{x}^\star_t; t) \geq 0$ for all $\bs{y} \in \Omega$, and therefore, since $\bs{x}^\star_{t-1} \in \Omega$ as it is a \gls{v-GNE} at $t - 1$, $\Delta_t^{\star^\top} F(\bs{x}^\star_t; t) \leq 0$. Thus, by adding and subtracting the term $\Delta_t^{\star^\top} G(\bs{x}^\star_{t-1}; t-1)$, we obtain
$
\Delta_t^{\star^\top} G(\bs{x}^\star_{t-1}; t-1) \leq \Delta_t^{\star^\top} (G(\bs{x}^\star_{t-1}; t-1) \!-\! (G(\bs{x}^\star_t; t) + U(\bs{x}^\star_t; t))).
$
With the personalized functionals in \eqref{eq:pers_feedback} and the perfect estimate of the pseudo-gradients, we have that $U(\bs{x}^\star_t; t) = c(t) (\Delta_t^\star + \alpha(t) G(\bs{x}^\star_{t-1}; t-1))$.
Then, in view of Standing Assumption~\ref{standing:symmetry}, it follows that 
$
\alpha(t) \Delta_t^{\star^\top} \nabla \theta(\bs{x}^\star_{t-1}; t -1) \leq \Delta_t^{\star^\top} ( \nabla \theta(\bs{x}^\star_{t-1}; t-1) + c(t) \bs{x}^\star_{t-1}  - \nabla \theta(\bs{x}^\star_t; t) - \xi(t) \bs{x}^\star_{t}).
$
Now, by summing and subtracting the term $\Delta_t^{\star^\top} \nabla \theta(\bs{x}^\star_{t}; t -1)$ in the RHS, we obtain
$$
\begin{aligned}
	&\alpha(t) \Delta_t^{\star^\top} \nabla \theta(\bs{x}^\star_{t-1}; t -1) \leq  \Delta_t^{\star^\top}(\nabla \theta(\bs{x}^\star_{t}; t) - \nabla \theta(\bs{x}^\star_{t}; t -1))\\
	&\hspace{5cm}- \Delta_t^{\star^\top} ( \nabla \theta(\bs{x}^\star_{t-1}; t-1) \!+\! c(t) \bs{x}^\star_{t-1}  \!-\! \nabla \theta(\bs{x}^\star_t; t-1) \!-\! c(t) \bs{x}^\star_{t}).
\end{aligned}
$$
As in the proof of Lemma~\ref{lemma:desc_perfect}, the second term in the RHS is always upper bounded by $-\ell \| \Delta_t^\star \|^2$ once introduced the auxiliary function $\psi(\bs{x};t) \coloneqq \theta(\bs{x}; t) + \tfrac{c(t)}{2} \| \bs{x} \|^2$ and chosen $c(t)$ in such a way that $\bs{x} \mapsto \psi(\bs{x};t)$ is $\ell$-strongly convex, i.e., $c(t) \geq 2 \ell$ for all $t \in \N$, in view of the $\ell$-weak convexity of $\bs{x} \mapsto \theta(\bs{x}; t)$. On the other hand, the first term attains its maximum positive module when the two vectors $\Delta_t^\star$ and $(\nabla \theta(\bs{x}^\star_{t}; t) - \nabla \theta(\bs{x}^\star_{t}; t -1))$ are aligned, thus leading to the following relation implied by Lemma~\ref{lemma:grad_tv}:
\begin{align}\label{eq:descent_tv}
	\Delta_t^{\star^\top} \nabla \theta(\bs{x}^\star_{t-1}; t -1) & \leq -\tfrac{\ell}{\alpha(t)} \| \Delta_t^\star \|^2 + \tfrac{1}{\alpha(t)} \|\Delta_t^\star\| \mathrm{e}_{\nabla},
\end{align}
where we have imposed that $\alpha(t) = 1 - c(t) \xi(t) > 0$, i.e., $\xi(t) < 1/c(t)$.
Finally, adding and subtracting the term $(1/4 \alpha(t) \ell) \mathrm{e}^2_{\nabla}$ allows us to complete the square in the RHS, thus yielding to \eqref{eq:descent_ball_tv}, and hence concluding the proof.
\hfill\qedsymbol

\smallskip

\textit{Proof of Theorem~\ref{th:conv_tv}}: We start by combining the descent lemma \cite[Prop.~A.24]{bertsekas1997nonlinear}, Lemma~\ref{lemma:grad_tv} and \ref{lemma:desc_perfect_tv} (specifically, the relation in \eqref{eq:descent_tv}). For all $t \in \N$, we obtain
\begin{align}\label{eq:tv_descent_lemma}
	\theta(\bs{x}^\star_t;t\!-\!1) &\leq \theta(\bs{x}^\star_{t - 1}; t\!-\!1) \!+\! (\Delta_t^\star)^\top \nabla \theta(\bs{x}^\star_{t - 1}; t\!-\!1) \!+\! \tfrac{\ell}{2} \|\Delta_t^\star\|^2\\
	&\leq \theta(\bs{x}^\star_{t - 1}; t-1)   - \ell \tfrac{2 - \alpha(t)}{2 \alpha(t)} \|\Delta_t^\star\|^2 + \tfrac{1}{\alpha(t)} \|\Delta_t^\star\| \mathrm{e}_{\nabla}.\nonumber
\end{align}

Adding and subtracting the term $\theta(\bs{x}^\star_t;t)$ yields to 
$$
\begin{aligned}
	\theta(\bs{x}^\star_t;t) &\leq \theta(\bs{x}^\star_{t - 1}; t-1) + \theta(\bs{x}^\star_t;t) - \theta(\bs{x}^\star_t;t-1)- \ell \tfrac{2 - \alpha(t)}{2 \alpha(t)} \|\Delta_t^\star\|^2 + \tfrac{1}{\alpha(t)} \|\Delta_t^\star\| \mathrm{e}_{\nabla}\\
	&\leq \theta(\bs{x}^\star_{t - 1}; t-1) + |\theta(\bs{x}^\star_t;t) - \theta(\bs{x}^\star_t;t-1)|- \ell \tfrac{2 - \alpha(t)}{2 \alpha(t)} \|\bs{x}^\star_t - \bs{x}^\star_{t - 1}\|^2 + \tfrac{1}{\alpha(t)} \|\bs{x}^\star_t - \bs{x}^\star_{t - 1}\| \mathrm{e}_{\nabla}\\
	&\leq \theta(\bs{x}^\star_{t - 1}; t-1)  \!-\! \ell \tfrac{2 - \alpha(t)}{2 \alpha(t)} \|\Delta_t^\star\|^2 \!+\! \tfrac{1}{\alpha(t)} \|\Delta_t^\star\| \mathrm{e}_{\nabla} + \mathrm{e}_{\theta},
\end{aligned}
$$
where the latter inequality follows in view of Assumption~\ref{ass:tv} i). By defining $\beta(t) \coloneqq \ell (2 - \alpha(t))/2 \alpha(t)$, which is strictly positive if $\xi(t) \in [0,1)$ for all $t \in \N$, and focusing on the second and third terms in the RHS of this latter inequality, we can complete the square by adding and subtracting $(4 \alpha^2(t) \beta(t))^{-1}\mathrm{e}^2_{\nabla}$.  After few algebraic manipulations, we obtain the following inequality
\begin{equation}\label{eq:sample_complexity_structure_tv}
	\begin{aligned}
		\theta(\bs{x}^\star_t;t) \leq  \theta(\bs{x}^\star_{t - 1}; t-1)  &-\beta(t) \left(	\|\Delta_t^\star\| - \tfrac{1}{2 \alpha(t) \beta(t)} \mathrm{e}_{\nabla} \right)^2+ \tfrac{1}{4 \alpha^2(t) \beta(t)} \mathrm{e}^2_{\nabla} + \mathrm{e}_{\theta}.
	\end{aligned}
\end{equation}

From now on, we follow the proofline introduced in \cite[Th.~3.1]{davis2019stochastic}, \cite[Th.~1]{mai2020convergence}, as in \eqref{eq:sample_complexity_structure} we have obtained a typical structure characterizing the sample complexity analysis for sample-based stochastic algorithms (see, for instance, \cite[Eq.(1.3)]{davis2019stochastic}). 
Specifically, by fixing some $T \in \N$ and summing up the inequality above over $t \in \mc{T} \coloneqq \{1, \ldots, T\}$, 
for any global minimizer $\bs{x}^\star(T) \in \Theta(T)$, we obtain	
$$
\begin{aligned}
	\theta(\bs{x}^\star(T); T) \leq \theta(\bs{x}^\star_T; T) \leq \theta(\bs{x}^\star_{0}; 0) &+ \mathrm{e}^2_{\nabla} \sum_{t \in \mc{T}} \tfrac{1}{4 \alpha^2(t) \beta(t)} + T \mathrm{e}_{\theta}- \sum_{t \in \mc{T}} \beta(t) \left(	\|\Delta_t^\star\| - \tfrac{1}{2 \alpha(t) \beta(t)} \mathrm{e}_{\nabla} \right)^2.
\end{aligned}
$$

After moving $\theta(\bs{x}^\star(T); T)$ in the RHS and $\sum_{t \in \mc{T}} \beta(t) (	\|\Delta_t^\star\| - \tfrac{1}{2 \alpha(t) \beta(t)} \mathrm{e}_{\nabla} )^2$ in the LHS, we obtain:
$
\sum_{t \in \mc{T}} \beta(t) (	\|\Delta_t^\star\| - \tfrac{1}{2 \alpha(t) \beta(t)} \mathrm{e}_{\nabla} )^2 \leq |\theta(\bs{x}^\star_{0}; 0) - \theta(\bs{x}^\star(T); T)|  + \mathrm{e}^2_{\nabla} \sum_{t \in \mc{T}} \tfrac{1}{4 \alpha^2(t) \beta(t)} + T \mathrm{e}_{\theta}.
$
Moreover, to upper bound the term $|\theta(\bs{x}^\star_{0}; 0) - \theta(\bs{x}^\star(T); T)|$, which explicitly depends on the considered time horizon $T$, we exploit the following chain of inequalities:	
$$
\begin{aligned}
	|\theta(\bs{x}^\star_{0}; 0) - \theta(\bs{x}^\star(T); T)| \leq &|\theta(\bs{x}^\star_{0}; 0) \!-\! \theta(\bs{x}^\star(0); 0)| \!+\! \sum_{t \in \mc{T}} |\theta(\bs{x}^\star(t); t) \!-\! \theta(\bs{x}^\star(t); t-1)|\\
	&\hspace{3.35cm} + \sum_{t \in \mc{T}} |\theta(\bs{x}^\star(t); t-1) - \theta(\bs{x}^\star(t-1); t-1)|,
\end{aligned}
$$
where the inequality is obtained after adding and subtracting the terms $\sum_{t \in \mc{T}} (\theta(\bs{x}^\star(t); t-1) + \theta(\bs{x}^\star(t-1); t-1))$ in the absolute value, and upper bounding each resulting summation. In view of Assumption~\ref{ass:tv} i), the first summation in the RHS of the this latter inequality is upper bound by $T \mathrm{e}_{\theta}$, while the second one by $(\ell/2) T \mathrm{e}^2_{\delta}$. This follows by combining i) the $\ell$-weak convexity of the mapping $\bs{x} \mapsto \theta(\bs{x}; t-1)$ (Lemma~\ref{lemma:weak_conv}), ii) \cite[Lemma~2.1]{davis2019stochastic} with the optimality condition $\nabla \theta(\bs{x}^\star(t-1); t-1) = 0$, and iii) the bound postulated in Assumption~\ref{ass:tv} iii). Then, by defining $\Delta_{0} \coloneqq |\theta(\bs{x}^\star_{0}; 0) - \theta(\bs{x}^\star(0); 0)|$, which does not depend on the length of the horizon $T$, we obtain, with $\phi \coloneqq 2\mathrm{e}_{\theta} + \tfrac{\ell}{2} \mathrm{e}^2_{\delta}$,
\begin{equation}\label{eq:11}
		\sum_{t \in \mc{T}} \beta(t) \left(	\|\Delta_t^\star\| \!-\! \tfrac{1}{2 \alpha(t) \beta(t)} \mathrm{e}_{\nabla} \right)^2 \leq \Delta_0 \!+\! T \phi + \mathrm{e}^2_{\nabla} \sum_{t \in \mc{T}} \tfrac{1}{4 \alpha^2(t) \beta(t)}.
\end{equation}	
	
Since $\beta(t) >0$ for all $t \in \mc{T}$, we can leverage the Jensen's inequality \cite[Th.~3.4]{rockafellar1970convex} on the convex function $(\cdot)^2$ to lower bound the summation in the LHS. 
To this end, we normalize the coefficient by multiplying and dividing by $\bar{\beta} \coloneqq \sum_{t \in \mc{T}} \beta(t)$, and we define $\hat{\beta}(t) \coloneqq \beta(t)/ \bar{\beta}$, thus obtaining 	 
$
\bar{\beta} (\sum_{t \in \mc{T}} \hat{\beta}(t) (	\|\Delta_t^\star\| - \tfrac{1}{2 \alpha(t) \beta(t)} \mathrm{e}_{\nabla} )^2 ) \geq
\bar{\beta} (\sum_{t \in \mc{T}} \hat{\beta}(t)	(\|\Delta_t^\star\| - \tfrac{1}{2 \alpha(t) \beta(t)} \mathrm{e}_{\nabla} ) )^2.
$ 
Then, replacing this latter inequality into \eqref{eq:11}, 	 
$
	\bar{\beta} (\sum_{t \in \mc{T}} \hat{\beta}(t) (\|\Delta_t^\star\| - \tfrac{1}{2 \alpha(t) \beta(t)} \mathrm{e}_{\nabla} ) )^2 = \tfrac{1}{\bar{\beta}} (\sum_{t \in \mc{T}} {\beta}(t)	(\|\Delta_t^\star\| - \tfrac{1}{2 \alpha(t) \beta(t)} \mathrm{e}_{\nabla} ) )^2
	\leq \Delta_0 + T \phi  + \mathrm{e}^2_{\nabla} \sum_{t \in \mc{T}} \tfrac{1}{4 \alpha^2(t) \beta(t)}, 
$
and performing few algebraic manipulations directly yields to,		
$$
\begin{aligned}
	&\sum_{t \in \mc{T}} \beta(t) \left(\|\Delta_t^\star\| - \tfrac{1}{2 \alpha(t) \beta(t)} \mathrm{e}_{\nabla} \right) \leq \sqrt{\sum_{t \in \mc{T}} \left( \tfrac{\bar{\beta}}{4 \alpha^2(t) \beta(t)} \mathrm{e}^2_{\nabla} +\beta(t) \left(\Delta_0 + T \phi \right) \right)}.
\end{aligned}
$$	

By bringing the term $- \mathrm{e}_{\nabla} \sum_{t \in \mc{T}} (1/2 \alpha(t))$ in the RHS, note that the obtained inequality still holds true if we premultiply both sides by $1/T$. Additionally, we have $(1/T) \sum_{t \in \mc{T}} \beta(t) \|\Delta^\star_t\| \geq (1/T) \underline{\beta} \sum_{t \in \mc{T}}  \|\Delta^\star_t\|$, with $\underline{\beta} \coloneqq \textrm{min}_{t \in \mc{T}} \beta(t)$. This directly yields to the relation in \eqref{eq:seq_tv}.
\hfill\qedsymbol

\smallskip

\textit{Proof of Lemma~\ref{lemma:desc_imperfect_recon_tv}}: 	By replicating the initial steps as in the proof of Lemma~\ref{lemma:desc_perfect_tv}, due to Assumption~\ref{ass:recons_error_tv} we have, for all $t \geq \bar{t}$,
$$
\begin{aligned}
	&\alpha(t) \Delta_t^{\star^\top} \nabla \theta(\bs{x}^\star_{t-1}; t -1) \leq  c(t) \xi(t) \Delta_t^{\star^\top} \epsilon_{t-1}- \Delta_t^{\star^\top} ( \nabla \theta(\bs{x}^\star_{t-1}; t-1) + c(t) \bs{x}^\star_{t-1}  - \nabla \theta(\bs{x}^\star_t; t) - c(t) \bs{x}^\star_{t}).
\end{aligned}
$$

Summing and subtracting the term $\Delta_t^{\star^\top} \nabla \theta(\bs{x}^\star_{t}; t -1)$ in the RHS yields to
$
	\alpha(t) \Delta_t^{\star^\top} \nabla \theta(\bs{x}^\star_{t-1}; t -1) \leq  \Delta_t^{\star^\top}(\nabla \theta(\bs{x}^\star_{t}; t) - \nabla \theta(\bs{x}^\star_{t}; t -1))+ c(t) \xi(t) \Delta_t^{\star^\top} \epsilon_{t-1} - \Delta_t^{\star^\top} ( \nabla \theta(\bs{x}^\star_{t-1}; t-1) + c(t) \bs{x}^\star_{t-1}  - \nabla \theta(\bs{x}^\star_t; t-1) - c(t) \bs{x}^\star_{t}).
$
In case $c(t) \geq 2 \ell > 0$ for all $t \in \N$, the third term in the RHS is always upper bound by $- \ell \| \Delta_t^\star \|^2$ in view of the $\ell$-strongly convexity of the time-varying auxiliary function $\theta(\bs{x}; t) + \tfrac{c(t)}{2} \| \bs{x} \|^2$. Moreover, the first (resp., second) term attains its maximum positive module when the two vectors $\Delta_t^\star$ and $(\nabla \theta(\bs{x}^\star_{t}; t) - \nabla \theta(\bs{x}^\star_{t}; t -1))$ ($\Delta_t^\star$ and $\epsilon_{t-1}$) are aligned, thus leading to the following relation
$$
\begin{aligned}
	\alpha(t) \Delta_t^{\star^\top} \nabla \theta(\bs{x}^\star_{t-1}) &\leq  -\ell \| \Delta_t^\star \|^2 +  c(t) \xi(t) \|\Delta_t^\star\|\|\epsilon_{t-1}\|+ \|\Delta_t^\star\|\|\nabla \theta(\bs{x}^\star_{t}; t) - \nabla \theta(\bs{x}^\star_{t}; t -1)\| ,
\end{aligned}
$$
once assumed that $\xi(t) \geq 0$.
Following Lemma~\ref{lemma:desc_perfect_tv}, the second term is upper bound by $\|\Delta_t^\star\| \mathrm{e}_{\nabla}$, while, in view of Assumption~\ref{ass:recons_error_tv}, for any $t \geq \bar{t}$ and $K \geq K_1$, it holds that $c(t) \xi(t) \|\epsilon_{t-1}\| \leq (1- \alpha(t)) \mathrm{e}(K-1)$ with probability $1 - \delta_1$.
For all $t \geq \bar{t}$, imposing $\alpha(t) > 0$, i.e., $\xi(t) < 1/c(t)$, yields to
\begin{equation}\label{eq:descent_ball_imp_tv}
	\begin{aligned}
		\Delta^{\star^\top}_t  \nabla \theta(\bs{x}^\star_{t-1}; t -1) &\leq  - \tfrac{\ell}{\alpha(t)} \| \Delta_t^\star \|^2+ \tfrac{1}{\alpha(t)} \|\Delta_t^\star\| (\mathrm{e}_{\nabla} + (1-\alpha(t))\mathrm{e}(K-1)),
	\end{aligned}
\end{equation} 
defining $\sigma(K,t) \coloneqq \mathrm{e}_{\nabla} + (1-\alpha(t))\mathrm{e}(K-1)$ and completing the square in the RHS by adding and subtracting $(1/4\alpha(t)\ell) \sigma^2(K,t)$ yields to \eqref{eq:descent_ball_imp_tv_}.
\hfill\qedsymbol

\smallskip

\textit{Proof of Theorem~\ref{th:conv_recons_tv}}: 	Also in this case, we combine the descent lemma \cite[Prop.~A.24]{bertsekas1997nonlinear} and \eqref{eq:descent_ball_imp_tv} from Lemma~\ref{lemma:desc_imperfect_recon_tv}. In view of Assumption~\ref{ass:recons_error_tv}, for all $t \geq \bar{t}$ and $K \geq K_1$, we obtain
$$
\begin{aligned}
	\theta(\bs{x}^\star_t;t-1) &\leq \theta(\bs{x}^\star_{t - 1}; t-1)   - \ell \tfrac{2 - \alpha(t)}{2 \alpha(t)} \|\Delta_t^\star\|^2+ \tfrac{1}{\alpha(t)} \|\Delta_t^\star\| (\mathrm{e}_{\nabla} + (1-\alpha(t))\mathrm{e}(K-1)),
\end{aligned}
$$
with probability $1 - \delta_1$. By introducing $\beta(t) \coloneqq \ell (2 - \alpha(t))/2 \alpha(t)$, as in the proof of Theorem~\ref{th:conv_tv} we first add and subtract the term $\theta(\bs{x}^\star_t;t)$ in the LHS to obtain
$
	\theta(\bs{x}^\star_t;t)\leq \theta(\bs{x}^\star_{t - 1}; t-1) + |\theta(\bs{x}^\star_t;t) - \theta(\bs{x}^\star_t;t-1)| - \beta(t) \|\Delta_t^\star\|^2 + \tfrac{1}{\alpha(t)} \|\Delta_t^\star\| (\mathrm{e}_{\nabla} + (1-\alpha(t))\mathrm{e}(K-1)) \leq \theta(\bs{x}^\star_{t - 1}; t-1)   - \beta(t) \|\Delta_t^\star\|^2 + \tfrac{1}{\alpha(t)} \|\Delta_t^\star\| (\mathrm{e}_{\nabla} + (1-\alpha(t))\mathrm{e}(K-1)) + \mathrm{e}_{\theta}.
$
By focusing on the second and third terms in the RHS of this latter inequality, we can complete the square by adding and subtracting $(1/4 \alpha^2(t) \beta(t)) (\mathrm{e}_{\nabla} + (1-\alpha(t))\mathrm{e}(K-1))^2$, which is always possible if $\xi(t) \in [0,1/c(t))$ for all $t \in \N$. Therefore, it follows that also $\beta(t) > 0$, for all $t \in \N$ (and hence for $t \geq \bar{t}$). After few algebraic manipulations, by introducing $\sigma(K,t) \coloneqq \mathrm{e}_{\nabla} + (1-\alpha(t))\mathrm{e}(K-1)$ we obtain what follows
\begin{equation}\label{eq:sample_complexity_structure_imp_tv}
	\begin{aligned}
		\theta(\bs{x}^\star_t;t) &\leq \theta(\bs{x}^\star_{t - 1}; t-1)  + \tfrac{1}{4 \alpha^2(t) \beta(t)} \sigma^2(K,t)  + \mathrm{e}_{\theta}- \beta(t) \left(	\|\Delta_t^\star\| - \tfrac{1}{2 \alpha(t) \beta(t)}  \sigma(K,t)  \right)^2.
	\end{aligned}
\end{equation}
From now on, we follow the proofline introduced in \cite[Th.~3.1]{davis2019stochastic}, \cite[Th.~1]{mai2020convergence}, as in \eqref{eq:sample_complexity_structure_imp_tv} we have obtained a typical structure characterizing the sample complexity analysis for sample-based stochastic algorithms (see, for instance, \cite[Eq.(1.3)]{davis2019stochastic}). Specifically, by fixing some $T \in \N$ and summing up the inequality above over $t \in \mc{T} \coloneqq \{\bar{t} + 1, \ldots, \bar{t} + T \}$, in view of the possibly sporadic communication between the agents and the central coordinator, let $q(t)$ be the subscript indicating the number of samples adopted by the learning procedure $\mathscr{L}$ at the $t$-th outer iteration. In particular, by assuming to initialize Algorithm~\ref{alg:two_layer} with $q(0) = K_1$ samples, $q(t) = q(t-1) + 1$ if, at the $t$-th iteration, the agents communicate $p_{t}$ to the coordinator, $q(t) = q(t-1)$ otherwise. Note that the existence of some $K \in \N$ is guaranteed in view of Assumption~\ref{ass:recons_error_tv}. Then, for any global minimizer $\bs{x}^\star(T) \in \Theta(T)$, with probability $1-\delta_1$ we obtain
$
	\theta(\bs{x}^\star(T); T) \leq \theta(\bs{x}^\star_T; T) \leq \theta(\bs{x}^\star_{\bar{t}}; \bar{t}) - \sum_{t \in \mc{T}} \beta(t) (	\|\Delta_t^\star\| - \tfrac{1}{2 \alpha(t) \beta(t)} \sigma(K,t)  )^2
	+ \sum_{t \in \mc{T}} \tfrac{1}{4 \alpha^2(t) \beta(t)} \sigma^2(K,t) + T \mathrm{e}_{\theta}.
$
Then, by moving $- \sum_{t \in \mc{T}} \beta(t) (	\|\Delta_t^\star\| - \tfrac{1}{2 \alpha(t) \beta(t)} \sigma(K,t))^2$ in the LHS and $\theta(\bs{x}^\star(T); T)$ in the RHS of this latter inequality, we obtain
$
	\sum_{t \in \mc{T}} \beta(t) (	\|\Delta_t^\star\| - \tfrac{1}{2 \alpha(t) \beta(t)} \sigma(K,t))^2 \leq \theta(\bs{x}^\star_{\bar{t}}; \bar{t}) - \theta(\bs{x}^\star(T); T) + \sum_{t \in \mc{T}} \tfrac{1}{4 \alpha^2(t) \beta(t)} \sigma^2(K,t) + T \mathrm{e}_{\theta} \leq |\theta(\bs{x}^\star_{\bar{t}}; \bar{t}) - \theta(\bs{x}^\star(T); T)| + \sum_{t \in \mc{T}} \tfrac{1}{4 \alpha^2(t) \beta(t)} \sigma^2(K,t) + T \mathrm{e}_{\theta}.
$
Moreover, to upper bound the term $|\theta(\bs{x}^\star_{\bar{t}}; \bar{t}) - \theta(\bs{x}^\star(T); T)|$, which explicitly depends on the horizon $T$, we exploit the following inequality obtained by adding and subtracting $\sum_{t \in \mc{T}} (\theta(\bs{x}^\star(t); t-1) + \theta(\bs{x}^\star(t-1); t-1))$ in the absolute value:
$$
\begin{aligned}
	|\theta(\bs{x}^\star_{\bar{t}}; \bar{t}) - \theta(\bs{x}^\star(T); T)| &\leq |\theta(\bs{x}^\star_{\bar{t}}; \bar{t}) \!-\! \theta(\bs{x}^\star(\bar{t}); \bar{t})| \!+\! \sum_{t \in \mc{T}} |\theta(\bs{x}^\star(t); t) \!-\! \theta(\bs{x}^\star(t); t-1)|\\
	&\hspace{4cm} + \sum_{t \in \mc{T}} |\theta(\bs{x}^\star(t); t-1) - \theta(\bs{x}^\star(t-1); t-1)|.
\end{aligned}
$$

In view of Assumption~\ref{ass:tv} i), the first summation in the RHS of the this latter inequality is upper bound by $T \mathrm{e}_{\theta}$, while the second one by $(\ell/2) T \mathrm{e}^2_{\delta}$. This follows by combining i) the $\ell$-weak convexity of the mapping $\bs{x} \mapsto \theta(\bs{x}; t-1)$ (Lemma~\ref{lemma:weak_conv}), ii) \cite[Lemma~2.1]{davis2019stochastic} with the optimality condition $\nabla \theta(\bs{x}^\star(t-1); t-1) = 0$, and iii) the bound postulated in Assumption~\ref{ass:tv} iii). Then, by defining $\Delta_{\bar{t}} \coloneqq |\theta(\bs{x}^\star_{\bar{t}}; \bar{t}) - \theta(\bs{x}^\star(\bar{t}); \bar{t})|$, which does not depend on the length of the horizon $T$, we obtain
\begin{equation}\label{eq:1tv}
	\begin{aligned}
		\sum_{t \in \mc{T}} \beta(t) &\left(	\|\Delta_t^\star\| - \tfrac{1}{2 \alpha(t) \beta(t)} \sigma(K,t)  \right)^2
		\leq \Delta_{\bar{t}} + T \phi+ \sum_{t \in \mc{T}} \tfrac{1}{4 \alpha^2(t) \beta(t)} \sigma^2(K,t),
	\end{aligned}
\end{equation}
with $\phi \coloneqq 2 \mathrm{e}_{\theta} + \tfrac{\ell}{2} \mathrm{e}^2_{\delta}$. Since $\beta(t) >0$ for all $t \in \mc{T}$, we can leverage the Jensen's inequality \cite[Th.~3.4]{rockafellar1970convex} on the convex function $(\cdot)^2$ to lower bound the summation in the LHS. 
To this end, we normalize the coefficient by multiplying and dividing by $\bar{\beta} \coloneqq \sum_{t \in \mc{T}} \beta(t)$, and we define $\hat{\beta}(t) \coloneqq \beta(t)/ \bar{\beta}$, thus obtaining 	 
$
	\bar{\beta} (\sum_{t \in \mc{T}} \hat{\beta}(t) (	\|\Delta_t^\star\| - \tfrac{1}{2 \alpha(t) \beta(t)} \sigma(K,t) )^2 ) \geq \bar{\beta} (\sum_{t \in \mc{T}} \hat{\beta}(t)	(\|\Delta_t^\star\| - \tfrac{1}{2 \alpha(t) \beta(t)} \sigma(K,t) ) )^2.
$
Then, replacing this latter inequality into \eqref{eq:1tv}, 	 
$$
\begin{aligned}
	\bar{\beta} \left(\sum_{t \in \mc{T}} \hat{\beta}(t)	\left(\|\Delta_t^\star\| - \tfrac{1}{2 \alpha(t) \beta(t)} \sigma(K,t) \right) \right)^2 &= \tfrac{1}{\bar{\beta}} \left(\sum_{t \in \mc{T}} \beta(t)	\left(\|\Delta_t^\star\| - \tfrac{1}{2 \alpha(t) \beta(t)} \sigma(K,t) \right) \right)^2\\
	&\leq \Delta_{\bar{t}} + T \phi + \sum_{t \in \mc{T}} \tfrac{1}{4 \alpha^2(t) \beta(t)} \sigma^2(K,t),
\end{aligned}
$$	 
and performing few algebraic manipulations directly yields to,		
$$
\begin{aligned}
	\sum_{t \in \mc{T}} \beta(t)	&\left(\|\Delta_t^\star\| - \tfrac{1}{2 \alpha(t) \beta(t)} \sigma(K,t) \right) \leq \sqrt{\sum_{t \in \mc{T}} \left( \beta(t) \left(\Delta_{\bar{t}} + T \phi\right) + \tfrac{\bar{\beta}}{4 \alpha^2(t) \beta(t)} \sigma^2(K,t)\right)}.
\end{aligned}
$$	
By bringing the term $- \sum_{t \in \mc{T}} (1/2\alpha(t)) \sigma(K,t)$ in the RHS, note that the obtained inequality still holds true if we premultiply both sides by $1/T$, thus obtaining the average over $T$ in the LHS. Moreover, we have that $(1/T) \sum_{t \in \mc{T}} \beta(t) \|\Delta_t^\star\| \geq (1/T) \underline{\beta} \sum_{t \in \mc{T}} \|\Delta_t^\star\|$, with $\underline{\beta} \coloneqq \textrm{min}_{t \in \mc{T}} \beta(t)$, which directly yields to the relation in \eqref{eq:seq_recons_tv} with probability $1-\delta_1$.
\hfill\qedsymbol

\subsection{On learning the Lipschitz constant $\ell$}\label{sec:learning_lipschitz}
The convergence results shown in the paper require one to know the constant of weak convexity $\ell \coloneqq \sum_{i \in \mc{I}} \ell_i$ characterizing the unknown (dis)satisfaction function $\theta$. Given the noncooperative nature of the problem considered, however, it seems reasonable that each Lipschitz constant $\ell_i$ represents a private information held by each agent. On the other hand, as long as the coordinator is endowed with a learning policy $\mathscr{L}$ to estimate $\nabla \theta$, one may wonder if it is also possible to learn the associated Lipschitz constant $\ell$. Therefore, given $K \in \N$ agents' feedback, in this appendix we assume to have available some estimate $\hat{\ell}_K > 0$ of $\ell$, and then we address this privacy issue in the stationary setting (the time-varying one is similar and it is hence omitted in the interest of the paper length).

\begin{lemma}\label{lemma:desc_perfect_recon_lipschitz}
	Let Assumption~\ref{ass:recons_error} hold true for some fixed $\delta_1 \in (0,1]$, $c(t) \geq \gamma \hat{\ell}_K$, for some design parameter $\gamma > 0$, chosen so that $\gamma\hat{\ell}_{K-1} - \ell > 0$, and $K \geq K_1$, and $\xi(t) \in [0, 1/c(t))$ for all $t \in \N$. Then, for all $t \geq \bar{t}$, with the personalized incentives in $\eqref{eq:pers_feedback}$, we have
	\begin{equation}\label{eq:descent_ball_lipschitz}
		\begin{aligned}
			\Delta_t^{\star^\top} \nabla \theta(\bs{x}^\star_{t-1}) &\leq  \tfrac{\alpha(t)\kappa^2(t)}{\gamma\hat{\ell}_{K-1} - \ell} \mathrm{e}^2(K-1)- \left(\sqrt{\tfrac{\gamma\hat{\ell}_{K-1} - \ell}{\alpha(t)}} \|\Delta_t^\star\| - \kappa(t) \sqrt{\tfrac{\alpha(t)}{\gamma\hat{\ell}_{K-1} - \ell}} \mathrm{e}(K-1) \right)^2,
		\end{aligned}
	\end{equation} 
	with probability $1-\delta_1$.
	\QEDB
\end{lemma}
\begin{proof}
	The first part of the proof replicates the one of Lemma~\ref{lemma:desc_perfect_recon}, albeit in this case we can not rely on the $\ell$-strongly convexity of the auxiliary function $\psi(\bs{x}; t) \coloneqq \theta(\bs{x}) + \tfrac{c(t)}{2} \| \bs{x} \|^2$, since the Lipschitz constant $\ell$ is assumed to be unknown. In view of \cite[Lemma~2.1]{davis2019stochastic}, we have that $\Delta_t^{\star^\top} ( \nabla \theta(\bs{x}^\star_{t-1}) + c(t) \bs{x}^\star_{t-1}  - \nabla \theta(\bs{x}^\star_t) - c(t) \bs{x}^\star_{t}) = \Delta_t^{\star^\top} ( \nabla \psi(\bs{x}^\star_{t-1};t)  - \nabla \psi(\bs{x}^\star_t;t)) \leq -(c(t) - \ell) \|\Delta_t^\star\|^2$. 
	Then, by letting $c(t) \geq \gamma \hat{\ell}_{K-1}$ for all $t \geq \bar{t}$, $K \geq K_1$, relying on Assumption~\ref{ass:recons_error} and following the same reasoning for the two vectors $\Delta_t^\star$ and $\epsilon_{t-1}$ as in the proof of Lemma~\ref{lemma:desc_perfect_recon}, directly yields 
	$$
		\Delta_t^{\star^\top} \nabla \theta(\bs{x}^\star_{t-1}) \leq - \tfrac{\gamma\hat{\ell}_{K-1} - \ell}{\alpha(t)} \| \Delta_t^\star \|^2 + \tfrac{1-\alpha(t)}{\alpha(t)} \| \Delta_t^\star \| \mathrm{e}(K-1),
	$$
	with probability at least $1 - \delta_1$. Adding and subtracting  $(\alpha(t)\kappa^2(t)/(\gamma \hat{\ell}_{K-1} - \ell)) \mathrm{e}^2(K-1)$ to complete the square in the RHS, with $\gamma$ such that $\gamma \hat{\ell}_{K-1} - \ell > 0$, leads to \eqref{eq:descent_ball_lipschitz}.
\end{proof}
Two considerations are now in order: i) $\mathrm{e}^2(K-1)$ prevents $\Delta_t^\star$ from being a descent direction for $\theta$, and ii) the condition $\gamma\hat{\ell}_{K-1} - \ell > 0$ requires one to implicitly overestimate the unknown  $\ell$. To overcome this latter issue, in the spirit of Assumption~\ref{ass:recons_error} we postulate the following condition:
\begin{assumption}\label{ass:bar_K}
	For any $\delta_2 \in (0,1]$, there exists some $K_2 < \infty$ such that $\prob\{\hat{\ell}_K / \ell \in [\rho_{\normaltext{\textrm{min}}}, \, \rho_{\normaltext{\textrm{max}}}] \mid \forall K \geq K_2\} \geq 1 - \delta_2$, for some $0 < \rho_{\normaltext{\textrm{min}}} \leq 1 \leq \rho_{\normaltext{\textrm{max}}}$. 
	\QEDB
\end{assumption}
With arbitrary probability $1 - \delta_2$, the ratio of $\hat{\ell}_K >0$ and $\ell > 0$ is thus upper and lower bounded by some known terms, $\rho_{\textrm{min}}$ and $\rho_{\textrm{max}}$. As Assumptions~\ref{ass:recons_error} and \ref{ass:recons_error_tv}, Assumption~\ref{ass:bar_K} is mild for a number of learning strategies, as LS~\cite{notarnicola2020distributed}. 

Without restrictions, we henceforward set $\delta_1 = \delta_2 = \delta/2$, thus requiring to satisfy Assumption~\ref{ass:recons_error} and \ref{ass:bar_K} with the same probability $1-\delta/2$, and their intersection with probability $1-\delta$ (via the union bound).

The following bound characterizes the sequence of \gls{v-GNE} originated by Algorithm~\ref{alg:two_layer} in case the Lipschitz constant of $\nabla \theta$, $\ell$, is not  available a-priori.
\begin{theorem}\label{th:stat_lip}
	Let Assumption~\ref{ass:recons_error}--\ref{ass:bar_K} hold true for some fixed $\delta \in (0,1]$, $\gamma \in (1/2\rho_{\normaltext{\textrm{min}}}, 1/\rho_{\normaltext{\textrm{min}}}]$, $c(t) \geq \gamma \hat{\ell}_{K-1}$ for some $K \geq \bar{K} \coloneqq \normaltext{\textrm{max}} \{K_1, K_2\}$, and $\xi(t) \in [2 (1 - \gamma \rho_{\normaltext{\textrm{min}}})/c(t), 1/c(t))$, for all $t \in \N$.
	Moreover, let some $T \in \N$ be fixed, $\mc{T} \coloneqq \{\bar{t}+1, \ldots, T + \bar{t}\}$ and, for any global minimizer $\bs{x}^\star \in \Theta$, $\Delta_{\bar{t}} \coloneqq \theta(\bs{x}^\star_{\bar{t}}) - \theta(\bs{x}^\star)$. 
	Then, with the personalized incentives in \eqref{eq:pers_feedback}, the sequence of \normaltext{\gls{v-GNE}} $(\bs{x}^\star_t)_{t \in \mc{T}}$, generated by Algorithm~\ref{alg:two_layer}, satisfies the following relation
	\begin{equation}\label{eq:seq_recons_lipschitz}
		\begin{aligned}
					\tfrac{1}{T} \sum_{t\in \mc{T}} \|	\Delta_t^\star	\| &\leq \tfrac{1}{T \underline{\beta}} \sum_{t \in \mc{T}}  \kappa(t) \mathrm{e}(q(t)) +\tfrac{1}{T \underline{\beta}} \sqrt{\sum_{t \in \mc{T}} \left( \beta_{q(t)}(t) \Delta_{\bar{t}} + \tfrac{\bar{\beta} \kappa^2(t)}{\beta_{q(t)}(t)} \mathrm{e}^2(q(t)) \right)}
		\end{aligned}
	\end{equation}
	with probability $1 - \delta$, where $\beta_{q(t)}(t) \coloneqq (2 \gamma \hat{\ell}_{q(t)} - \ell (2 + \alpha(t)))/2 \alpha(t)$, $\bar{\beta} \coloneqq \sum_{t \in \mc{T}} \beta_{q(t)}(t)$, $\underline{\beta} \coloneqq \normaltext{\textrm{min}}_{t \in \mc{T}} \beta_{q(t)}(t)$, and $q(t) \geq \bar{K}$ is the number of available agents' feedback at the $t$-th outer iteration, $t \in \mc{T}$.
	\QEDB
\end{theorem}
\begin{proof}
	The proof mimics the same steps of the one of Theorem~\ref{th:conv_recons}. However, to recover the traditional structure characterizing the sample complexity analysis for sample-based stochastic algorithms as in, e.g., \cite[Eq.(1.3)]{davis2019stochastic}, we shall guarantee first that $\beta_{K-1}(t) \coloneqq (2 \gamma \hat{\ell}_{K-1} - \ell (2 + \alpha(t)))/2 \alpha(t)$, which follows by combining the descent lemma and \eqref{eq:descent_ball_lipschitz}, is strictly positive, for all $t \geq \bar{t}$, $K \geq \bar{K} \coloneqq \textrm{max} \{K_1, K_2\}$. Therefore, $\beta_{K-1}(t) > 0$ if and only if $\alpha(t) < 2 (\gamma (\hat{\ell}_{K-1} / \ell) - 1)$ and $\xi(t) < 1/c(t)$. Note that, due to Assumption~\ref{ass:bar_K}, the former inequality is always satisfied if $\alpha(t) < 2 (\gamma \rho_{\textrm{min}} - 1)$, as $\hat{\ell}_{K-1} / \ell \geq \rho_{\textrm{min}}$ with probability $1 - \delta$ for all $K \geq \bar{K} \geq K_2$. Thus, in view of the definition of $\alpha(t)$, we have that $\xi(t) > 2(1 - \gamma \rho_{\textrm{min}}) / c(t)$, for all $t \geq \bar{t}$, which amounts to a nonnegative lower bound in case $\gamma \leq 1/\rho_{\textrm{min}}$, whereas it shall verify $\gamma > 1/2\rho_{\textrm{min}}$ to meet $\xi(t) < 1/c(t)$. From now on, the proof is a verbatim copy of the one of Theorem~\ref{th:conv_recons} with $\beta_{K-1}(t)$ instead of $\beta(t)$, while in summing up over $\mc{T}$ we account for the sporadic communication among agents and central coordinator by considering $\beta_{q(t)}(t)$.
\end{proof}
The considerations following Theorem~\ref{th:conv_recons}, i.e., the parameter tuning for the personalized functionals, big-$O$ analysis and learning strategy $\mathscr{L}$ apply also to this case \emph{mutatis mutandis}. 
\bibliographystyle{IEEEtran}
\bibliography{21_TAC.bib}



\end{document}